\documentclass{amsart}
\usepackage{amsthm,amsfonts,amsmath,amscd,amssymb,latexsym,epsfig,eucal,float}
\usepackage{multirow}

\newtheorem{theorem}{Theorem}
\newtheorem{lemma}[theorem]{Lemma}
\newtheorem{corollary}[theorem]{Corollary}
\newtheorem{proposition}[theorem]{Proposition}
\newenvironment{remark}{\medskip \refstepcounter{theorem}
\noindent  {\bf Remark \thetheorem}.\rm}{\,}

\def\BOne{{\mathchoice {\rm 1\mskip-4mu l} {\rm 1\mskip-4mu l}
                          {\rm 1\mskip-4.5mu l} {\rm 1\mskip-5mu l}}}

\def\<{\langle}

\def\>{\rangle}

\def\mb#1{{\mathbb #1}}
\def\mc#1{{\mathcal #1}}

\begin{document}

\title[Vibrations of thin plates]{The (small) vibrations
of thin plates \\ {\it \footnotesize In memory of Stephen A. Andrea}}
\author{Santiago R. Simanca}
\address{888 S Douglas Road, Apt 121, Coral Gables, FL 33124, U.S.A.}
\email{srsimanca@gmail.com}

\begin{abstract}
We describe the equations of motion of  
elastodynamic bounded bodies in 3-space, and their linearizations at a 
stationary point. Using the latter as an approximation to model  
small motions, we develop a scheme to find numerical solutions of these 
equations. We discretize the solution in the space of PL vector fields 
associated to the oriented faces of the first barycentric subdivision of a 
given smooth initial triangulation of the body, in order to exploit the 
algebraic topology properties of the body that these vector fields encode 
into the sought after solution, and solve a weak version of the linearized 
equations in that context. We apply our scheme to a couple of 
relevant examples of thin bodies, bodies where one of the dimensions is at 
least one order of magnitude in size less than the  other two, and determine 
numerical approximations to some of their resonance modes of vibration. The
results obtained are consistent with known vibration patterns for 
these bodies derived experimentally.
\end{abstract}

\maketitle 

\section{Introduction}
It is a classical experimental result to obtain the normal modes of vibrations 
of an unassembled violin plate by horizontally mounting the plate with its 
belly down, sprinkling the other side with a fine powder, and causing the 
plate to vibrate by applying an external sinusoidal pressure force that points 
in the vertical direction, and that hits the belly of the mounted plate 
straight up. At certain frequencies resonance occurs, and the vibration of 
the plate bounces 
the powder into the nonvibrating nodal areas, outlining the nodal and 
antinodal configurations of the plate at its eigenfrequency modes. This 
experiment goes as far back as F\'elix Savart in 1830, building on a method 
developed by Ernst Chladni, (see, for instance, \cite[$2^{\rm nd}$ 
column, page 171]{hu}; this reference describes various other types of acoustic 
experiments also).  Today, the nodal lines of the plate vibrating at these 
frequencies are called the Chladni lines, and the entire portrait 
that they make is called the Chladni pattern.    

In this article we study the theoretical underpinning of these experiments. We
consider the plate material as an incompressible orthotropic elastic 
Hooke body, whose stored energy function is characterized by nine independent 
elastic constant parameters. With the density as an additional parameter, we 
describe the equations of motion governing the time evolution of these bodies,
and concentrate our attention on solving these equations numerically. We then
use these solutions to describe the Chladni patterns of the vibrating violin 
plate, and more generally, the Chladni patterns of very thin plates.   

The Cauchy problem for the equations of motion of an incompressible elastic 
coercive body whose boundary is free to move is well-posed, at least for 
a short time \cite{ebsi,ebsi2}, and the derivation of numerical solutions to it
is a matter of interest in its own right. This is a 
difficult problem, essentially because the nonlinearity 
of the equation of motion being solved is of a pseudodifferential nature,
with changes to the solution in a neighborhood of any point affecting the
solution everywhere else, all at once, since the solution of the pressure 
equation depends on global conditions across the body.

We solve these difficulties in two steps: The first is to use discretizing 
spaces for the solution that we seek that have encoded into them the algebraic
topology of divergence-free fields. Such a condition,    
a matter of satisfying an equation that locally involves finitely many of the 
coefficients in the discretized unknown, is then left to be regulated by
the equation itself, and if holding initially, it  
should be satisfied at later times within a small margin of error also due to 
the well-posedness of the problem we solve, thus ensuring that the nonlocal 
incompressible condition holds in time. The discretizing spaces used are 
associated to the dual of the Whitney forms of an oriented triangulation of 
the body, which in dimension 
three leads to a simplicial complex with functions as complex groups in 
degrees zero, and three, and vector fields in degrees one, and two, 
respectively, whose cohomology is the cohomology of the body, 
and which in degree two have as cycles the Abelian group of divergence-free 
vector fields where the actual solutions to the equations of motion lie; the 
second step arises by observing that the Chladni patters must 
involve primarily the small vibrations of the bodies under consideration, and 
these can be approximated well by the solution to the linearized equations of
motion about the canonical stationary state. Under that assumption, which
we follow, we can tame the effects of the almost degenerate nature of 
the bodies, bounded three dimensional with one of its directions at least 
one order of magnitude smaller in length than the other two, by choosing 
the simplices in the triangulation to be nondegenerate, and using sufficiently
many of them so that the weak formulation of the problem that we solve in
our numerical scheme captures well the vibrations of the body at any point.
The solution so obtained feels the contributions to the vibration modes 
arising from all of the points, including the many far apart boundary points 
that are separated by a very small distance within the body. It does so, 
however, at a very large computational memory expense, a consequence of the 
need to triangulate the body with the appropriate resolution in order to 
derive accurate results. 

The large systems of equations that need to be solved, however sparse, are 
computationally challenging, and this is 
a difficulty that we must face even after resolving that of determining the 
systems themselves. The judicious choice of discretizing spaces minimizes 
these difficulties. Since we encode in them the initial conditions of the
waves we seek, with  the algebraic topological properties that they have,
we can then allow the well-posedness of the equations solved for to play its 
role in keeping the numerical solutions derived physically accurate in time.
  
The way the various types of waves propagate within the 
body is a property encoded in the density and elastic constants parameters
that are part of the equations of motion. We may think of the body as an
ensemble of springs, one per every pair of faces in the triangulation whose
stars intercept nontrivally, placed within the body according 
to the orientation of the said faces. The waves within the body, Lamb, 
Rayleigh, shear, or otherwise, result from the interaction of the elements
of this ensemble of springs as they move, which is ruled by the 
equations of motion. They all fit now within a single framework.
And the solutions we construct numerically to describe these waves are built 
to maintain the global topological constraint imposed by the 
incompressibility condition in time, yielding very accurate approximations to 
the actual motion of the body. The innovative use of our discretizing 
spaces counterbalances the need to go through the computational complexities 
inherent in the problem, in exchange for producing results that are faithful 
to the physical reality of the motion while in the elastic regime.

Although here we apply our approach to orthotropic bodies only, the same 
could be done to study the vibrations of elastic bodies in general, as long as
their stored energy function is (or is assumed to be) smooth. The variational 
principle used to derive the equations of motion can be extended to treat 
cases where part of the boundary is fixed, while the rest is free to move, and
so, for instance, we could treat the problem of 
isotropic functionally graded plates \cite{mahi,benn} in a way that is far 
easier by comparison with the problems analyzed here,
given the much smaller number of elastic constant parameters even if their 
values now vary across the body. The 
resulting simulations would be computationally more complex than approaches 
describing the motion in terms of an ad hoc small number of degrees of 
freedom, and basic assumptions, but the physical meaning of the results  
would be unquestioned.
 
The simulations we carry out fit well with the classical experiments on
violin plates, and serve also in the study of small vibrations of thin plates 
in general, as observed above. They can be used to solve problems of practical 
importance. For assuming given a set of parameters determining the elastodynamic
properties of the material being used in making the thin plate, and a desired 
Chladni pattern, we may find the geometry of the plate that conforms to that 
pattern. This gives an idea on how to carve the material in order 
to achieve a specified vibration pattern, or how this pattern is affected
when the elastodynamic constants of the material vary across it, 
for instance, the material that results by coating with a protective layer 
of varnish a given plate of wood. We can find the natural vibration pattern
of a sheet made of aluminum, and study how that will change when instead
we use an alloy of aluminum and copper to make the sheet with a material 
whose strength is larger. In fact, given the elastic constants of the 
material at room temperature, we should be able to determine the natural 
free and dampening vibration modes of a pointy thin cylindrical shell made 
out of this material when accelerating in 
open space, and improve on the result significantly if we take into 
consideration the rising temperature of the accelerating shell that is 
due to dampening modes, and consider its effect on the elastic constants of 
the material used to make it as time passes by, again producing results that
are based on the fundamental physical principles ruling the motion.  

\subsection{Organization of the article}
In \S\ref{s2} we state the equations of motion of incompressible 
elastodynamic bodies, and briefly sketch the proof of the well-posedness 
of its Cauchy problem in $C^{k,\alpha}$ spaces, spaces which are 
better suited for our purpose here than the Sobolev spaces used in the 
original proof \cite{ebsi,ebsi2}. We derive also the specific linearization 
of these equations that we wish to solve numerically here, with its 
natural set of boundary 
conditions. In \S\ref{s3}, we define the generalized Hooke materials, 
among which are the orthotropic ones that we use, and after recalling the 
general notions of stress and strain tensors. For the sake of completeness, 
we recall also the definition of the elastic parameters that characterize 
their stored energy functions, and the data for the material that we shall 
use in our 
numeric simulations in \S\ref{s5} below. In \S\ref{s4}, we write down the 
algorithms that we use to 
solve the said equation of \S\ref{s2}, after a careful discussion of 
the discretizing spaces, closely related to the divergence-free vector
fields that represent cohomology classes of the body in degree two, and
which we describe using the metric dual of the Whitney forms associated to 
the faces of the first barycentric subdivision of a given triangulation.
We emphasize the appropriateness of this choice by describing the metric 
dual simplicial complex of the usual de Rham complex, all of whose Abelian
groups of chains can be expressed in terms of the metric duals of the
Whitney forms of various degrees, and show how the divergence-free vector
fields appear 
as the group of cycles of this complex in degree two. 
In \S\ref{s5}, we present the various simulations that we have carried out,
as well as some additional numerical computations to support the virtues
of our approach. We end with some remarks in \S\ref{s6}, where among other 
things, we observe that the solution to the linearized equations of motion, 
which we study here as approximations to small vibrations of the body, is 
essentially the initial step in a Newton iteration scheme to solve 
the actual equations of motion, at least on the time interval where they
do not develop singularities (if at all). 

\section{The equations of motion}
\label{s2}
In this section, we begin by recalling the equations of motion of an 
elastodynamic body, and their linearization.
We refer the reader to \cite{ebsi,ebsi2} for details.

Let us suppose that $\Omega$ is a bounded domain in $\mb{R}^3$. 
We assume its boundary to be of class $C^{1,\alpha}$,
and that the density $\rho$ of the material filling $\Omega$ is 
constant.  

The motion of $\Omega$ is encoded into a curve
$$
\eta(t) : \Omega \hookrightarrow \mb{R}^3
$$
of embeddings of $\Omega $ into $\mb{R}^3$ that preserve 
volume. The path $t \rightarrow \eta (t)(x)$ denotes the position at 
time $t$ of a particle initially at $x\in \Omega$, and if  
$D \eta(t)(x)= (\partial_{x^i}\eta^j(t)(x))$ is the deformation gradient, 
the volume preserving property or {\it incompressibility}, is equivalent to
\begin{equation} \label{eq1}
J(\eta(t))=\det{D \eta (t)(x)}= 1\, .
\end{equation}

The material properties of $\Omega$ are characterized by its 
stored energy function $W$. This function indicates how the various 
particles of $\Omega$ are bound to each other, and it is often assumed
to be a function of the deformation gradient, $W=W(D\eta)$. A 
 typical choice is that given by $W(\partial _{i}\eta ^{j})={\displaystyle 
-\frac{3}{2}+ \frac{1}{2}}
\partial _{i}\eta ^{\alpha } \partial _{i}\eta ^{\alpha }$ 
which models isotropic neoHooke materials (throughout the paper   
we use the standard convention of adding over repeated indexes). 
For the time being, we keep a general point of view in mind.

In the presence of no external forces, the trajectory of the body
will be an extremal path of the Lagrangian
$$  
{\mc{L}}(\eta )=\frac{1}{2}\int _{0}^{T}\int _{\Omega }
\rho \| \dot{\eta } (t)(x)\| ^{2}dxdt - \int _{0}^{T}\int _{\Omega }
W(\partial _{i} \eta ^{\alpha }(t)(x))dxdt\; .
$$  
The stationary points of $\mc{L}(\eta)$ must be searched for among variations 
of $\eta(t)$ preserving (\ref{eq1}). If we let 
$\eta (t,s)$ be one such variation, and define  
$\zeta: \eta(t)(\Omega) \rightarrow \mb{R}^{3}$ by
$$
\zeta(\eta(t)(x))=\partial _{s} \eta (t,s)\mid _{s=0}(x)\, ,
$$
we then have that ${\rm div}\, \zeta =\partial_{x^j}\zeta^j=0$, and  
$$
\begin{array}{rcl}
{\displaystyle \frac{d}{ds}\mc{L}(\eta (t,s))\mid _{s=0}} & = &
{\displaystyle \int _{0}^{T}\int _{\Omega } \langle (-\rho\ddot{\eta }+{\rm Div}
W^{'}(D\eta ))(x), \zeta(\eta(t)(x) \> dxdt}
\\ &  & {\displaystyle \;\;\;-\int _{0}^{T}
\int _{\partial \Omega }\< W^{'}(D\eta )N,\zeta(\eta(t)(x))\> d\sigma dt } \; , 
\end{array} 
$$
where $W^{'}$ means the derivative of $W$ with respect to the 
variables $D\eta $, ${\rm Div}$ is the divergence of $W^{'}$ with respect to
the material coordinates $x$, $N$ is the unit normal to $\partial \Omega$ and
$d\sigma $ is the surface measure for $\partial \Omega $. It thus follows that
the motion of the body takes place along a solution of the free boundary value
problem
\begin{equation}\label{eq2} 
\rho \ddot{\eta }(t)(x)-{\rm Div}\, W^{'}(D\eta (t)(x))=
\nabla p(t)(\eta (t)(x)) \, ,
\end{equation}  
\begin{equation}
W^{'}(D\eta )N+p(t)(\eta(t)(x))J^{b}(\eta)\nu \circ \eta =0 \; {\rm on}
\; \partial \Omega \, , 
\label{eq3} 
\end{equation}
for a certain function  
$p(t): \eta (t)(\Omega ) \rightarrow \mb{R}$.
Here, $J^b(\eta(t))$ is the Jacobian determinant of $\eta(t)$ restricted to 
the boundary $\partial \Omega$, and $\nu: \partial (\eta(t)(\Omega)) 
\rightarrow \mb{R}^3$ is the unit normal. 

In coordinates $\eta^\alpha=\eta ^{\alpha }(x^{1}, x^2 , x^{3})$, we have that
$$ 
({\rm Div}\,  W^{'})^{\alpha }
=\partial _{x^{i}}\left( \frac{\partial W}{\partial (\partial_{i}\eta ^{\alpha
}
)}\right)\stackrel{{\rm def}}{=} A_{ij}^{\alpha \beta }(D\eta )\partial _{i}
\partial _{j}\eta ^{\beta }
$$ 
and
$$ 
(W^{'}(D\eta )N)^{\alpha }=
\frac{\partial W}{\partial (
\partial_{i}\eta ^{\alpha })}N^{i}\, , 
$$
respectively,  and the system above is given by 
$$
\rho \ddot{\eta }^{\alpha }  =  A_{ij}^{\alpha \beta }(D\eta )\partial _{i}
\partial _{j}\eta ^{\beta }+(\partial _{\alpha }(p(t)))\circ \eta \, , 
$$

$$
\frac{\partial W}{\partial (
\partial_{i}\eta ^{\alpha })}N^{i} +p(t)(\eta(t)(x))J^b(\eta) (\nu\circ \eta)^
\alpha  = 0 \, . 
$$

We require that the stored energy function $W$ be {\it coercive} at  
least for small vibrations. That is to say, 
we require the operator $A_{ij}^{\alpha \beta}(D\eta)\partial_i \partial_j$ to
be uniformly elliptic in a neighborhood of the stationary curve $\eta(t)(x)=x$.

Equations (\ref{eq1}), (\ref{eq2}), (\ref{eq3}) are the equations of motion 
of an incompressible elastic body of stored energy $W$.

Given an operator $F$, we define the operator $F_{\eta}$ by 
$F_{\eta}u=(F(u\circ \eta^{-1})) \circ \eta$. If we now differentiate 
(\ref{eq1}) with respect to $t$, we obtain that 
${\rm div}_{\eta}\dot{\eta}=0$, and upon a second differentiation, we have that
$$
{\rm div}_{\eta} \ddot{\eta}= -[ (\dot{\eta}\circ \eta^{-1}\cdot \nabla, {\rm div}]_{\eta}
\dot{\eta}={\rm trace} (D_{\eta}\dot{\eta})^2 \, .
$$
We use this result to express equations (\ref{eq1}), (\ref{eq2}), (\ref{eq3}) 
as a first order system.

\begin{proposition}
Let $\eta(t)$ be a curves of embeddings describing the motions of an 
incompressible elastic
body of constant density $\rho$ and stored energy function $W$. Then 
${\rm div}_{\eta} \dot{\eta}=0$, the vector $W'(D\eta)N$ is 
perpendicular to $W'(D\eta)T$ {\rm (}where $T$ is any vector tangent to $\partial \Omega${\rm )},
and 
\begin{equation}
\frac{d}{dt}\left( \begin{array}{c}
                   \eta \\ \rho \dot{\eta } 
                   \end{array}\right) = \left(
                                            \begin{array}{c}
                                            \dot{\eta } \\ A(\eta ,\dot{\eta })
                                        \end{array}\right) 
        \stackrel{{\rm def}}{=} F(\eta ,\dot{\eta }) \; ,
\label{eqs} 
\end{equation}
where
$$
A(\eta ,\dot{\eta }) = {\rm Div}\, W'(D\eta) + \nabla _{\eta }q \, ,
$$ 
and $q$ solves the boundary value problem
\begin{equation}
\begin{array}{c@{=}l}
\Delta _{\eta }q \; & \; -{\rm div}_{\eta }{\rm Div}\, W'(D\eta) +\rho \, 
{\rm trace}(D_{\eta } \dot{\eta })^{2} \vspace{1mm} \\
q\mid _{\partial \Omega }\; & \;-{\displaystyle \frac{\langle W'(D\eta) N
 , \nu \circ \eta  \rangle }{J^{b}(\eta )}} \; .
\end{array}
\label{eqb} 
\end{equation}
Here, $q=p\circ \eta $ where $p$ is as in $(\ref{eq2})$, $N$ and $\nu$ are the
unit vectors normal to $\partial \Omega$ and $\partial (\eta(t)(\Omega))$, 
respectively, and $J^b(\eta)$ is the Jacobian determinant of $\eta$ 
restricted to the boundary. 
\end{proposition}

The function $q$ is normally called the {\it pressure}.   

Since $W$ is coercive, the boundary value problem (\ref{eqb}) is elliptic, 
and has a unique pressure function solution $q=q(\eta,\dot{\eta})$. This 
solution is a nonlocal pseudodifferential operator in $(\eta, \dot{\eta})$.
The system (\ref{eqs}) admits the time independent curve 
$\eta=\BOne$ as a solution, with pressure function the constant $q=-\< W'(\BOne)N,N\>$.  

For reasonably smooth $\Omega$ and coercive stored energy function $W$,
the Cauchy problem for (\ref{eqs}) with initial condition 
$(\eta (0),\dot{\eta}(0))=(\BOne, w)$, $w$ a divergence-free vector field, 
is well-posed. This was proven in 
\cite{ebsi2} working on Sobolev spaces. With some additional difficulties, 
the same argument employed there  
can be cast working over $C^{k,\alpha}$-spaces instead; we just need to 
invoke the continuity of pseudodifferential operators on these spaces, as 
opposed to that on Sobolev spaces, and redo the said argument in this new 
context. The benefit of doing so is twofold: on the one hand, 
$C^{k,\alpha}$-spaces are better suited to analyze the question of 
consistency of the numerical solutions of the equations we derive in this 
article, and for which we use the PL vector fields associated to the 
elements of a triangulation as discretizing spaces; on the 
other hand, we get optimal regularity results for the analysis of the 
Cauchy problem, which can be started if we assume merely that $\eta(t)$ is 
a $C^{2,\alpha}$ curve of embeddings. It is for this reason that we
discuss the continuity properties of $F(\eta, \dot{\eta})$ over these other 
spaces, as opposed to those employed in \cite{ebsi,ebsi2}.

Let $C^{k,\alpha}(\Omega; \mb{R}^3)$ be the space of vector fields with 
$C^{k,\alpha}$-regularity. These spaces are closed under multiplication by 
scalar valued functions in $C^{k,\alpha}$, and pseudodifferential operators
of order $n$ act continuously from $C^{k,\alpha}$ into $C^{k-n,\alpha}$ 
\cite{tay}.  We define $\mc{C}_{k,\alpha}=C^{k,\alpha}(\Omega ; 
\mb{R}^{3})\oplus C^{k-1,\alpha}(\Omega ; \mb{R}^{3})$, and view 
$F$ in (\ref{eqs}) as a map:
\begin{equation}
 \begin{array}{cccl}
F: &  \mc{C}_{k+1,\alpha} & \longrightarrow &  \mc{C}_{k,\alpha} \\  & 
(\eta ,\dot{\eta })& \longrightarrow & F(\eta ,\dot{\eta })
\end{array}\; . 
\label{eq6} 
\end{equation}
In order to preserve the validity of the variational principle producing the 
equations of motion as a duality pairing, it is natural to impose the 
restriction that $k\geq 0$. We see below that the stricter assumption 
$k\geq 1$ is required in order to make classic sense of other terms. We 
shall assume the latter always.

The continuity of $F$ is a question of continuity of the map $A$ from 
$\mc{C}_{k+1,\alpha}$ to 
$C^{k-1,\alpha}(\Omega ; \mb{R}^{3})$. Thus, if $(\eta ,\dot{\eta })\in 
\mc{C}_{k+1,\alpha}$, 
the pressure function $q$  obtained by solving (\ref{eqb}) must be in 
$C^{k,\alpha}(\Omega )$.
If $\eta $ is a $C^{k+1,\alpha}$-diffeomorphism from $\Omega $ to $\eta 
(\Omega)$, $\dot{\eta }\in C^{k,\alpha}(\Omega; \mb{R}^{3})$, as $k\geq 1$, 
we have that
$-{\rm div}_{\eta }{\rm Div}\, W'(D\eta)+{\rm trace}(D_{\eta }\dot{\eta })^{2}$
is a well-defined distribution in $C^{k-2,\alpha}(\Omega )$. Since the 
boundary condition for
$q$ lies in $C^{k,\alpha}(\partial \Omega)$, the regularity of solutions
to the Dirichlet problem implies that $q$ is in $C^{k,\alpha}$, as desired.

\begin{theorem} {\rm (\cite[Theorem 5.53 and \S 6]{ebsi2})} \label{th2} 
Let $\Omega$ be an incompressible elastic body of constant density $\rho$ and 
coercive stored energy function $W$. Given an initial condition $(\eta(0),
\dot{\eta}(0)) =(\BOne ,w)$, $w\in C^{k,\alpha}(\Omega; \mb{R}^3)$, $k\geq 1$, 
${\rm div}\, w=0$, there exists a positive real number $T$ such that the 
initial value problem for {\rm (\ref{eq1}), (\ref{eq2}), (\ref{eq3})}
has a unique solution $\eta(t)\in C^2([0,T]; C^{k+1,\alpha}(\Omega,\mb{R}^3))$ 
with the said initial condition. The value of $T$ depends upon the norm 
of finitely many derivatives of $w$ only, and the solution $\eta(t)$ is a 
continuous function of the initial data.
\end{theorem}

Our interest is in the motion of the bodies when they experience relatively 
small deformations from the stationary solution $\eta(t)(x)=x$. We shall study 
numerically the time evolution  
of the solutions to the linearized equations about $\eta(t)(x)=x$, and use 
them as actual approximations to the nonlinear small vibrating motions.   
We end this section by outlining the important details of the derivation of 
this linearization. (The derivation in \cite[Proposition 2.8]{ebsi} was 
carried out for neoHooke materials only.) We linearize first about a 
fixed but arbitrary element $(\eta , \dot{\eta })$ of $\mc{C}_{k+1,\alpha}$ 
whose first component is a diffeomorphism onto its image, and then apply the
result at the said stationary point.  

We begin by linearizing the boundary condition (\ref{eq3}), which is      
required in computing the linearization of the operator $F$ in (\ref{eqs}).
We fix $x\in \partial \Omega $, and pick $\{ T_1, T_2\}$
orthonormal and tangent to $\partial \Omega $ at $x$. We fix also the 
coordinates in $\mb{R}^{3}$ so that $\partial_{T_{i}}\eta \in \mb{R}^2
\subseteq \mb{R}^{3}$ for $i=1,2$. Then,
$$
(J^{b}(\eta) \nu \circ \eta)(x)=\det{\left(
 \begin{array}{c}
\partial_{T_{1}}\eta \\ \partial_{T_{2}}\eta \end{array} \right)}
(0,0,1) \, . 
$$ 
We let $\eta(s)=\eta +su$ and $\dot{\eta}(s)=\dot{\eta}+sv$ be the variations of $\eta$ and 
$\dot{\eta}$ in the directions of $u$ and $v$, respectively. Then,
\begin{equation} \label{bd}
\displaystyle 
\< \partial _{s}(J^{b}(\eta) \nu \circ \eta), \nu \circ
\eta \>\mid_{s=0} = \\ J^{b}(\eta)\left[ \det{
      \left(
      \begin{array}{c}
      \partial_{T_{1}} u \\ \partial_{T_{2}}\eta 
      \end{array} 
      \right)}+
 \det{\left(  
     \begin{array}{c} 
     \partial_{T_{1}}\eta \\ \partial_{T_{2}} u 
     \end{array} 
      \right)}\right] \, .
\end{equation}
This is a function on $\partial \Omega $. Notice that the vector  
$\partial _{s}(J^{b}(\eta) \nu \circ \eta)$ itself may have a nontrivial
component tangent to $\eta (\partial \Omega )$. 

The boundary condition (\ref{eq3}) can be 
rewritten as $q(t)(x)J^b(\eta)=-\< W'(D\eta)N,\nu\circ \eta\>$.
If we replace $\eta$ by $\eta(s)$ in (\ref{eq3}), and differentiate with 
respect to $s$, we find that
$$
A_{ij}^{\alpha \beta }(D\eta )\partial _{j}u^{\beta }N^{i}\nu ^{\alpha }\circ
\eta +(\partial _{s}q)J^{b}(\eta) +q\< \partial _{s}(J^{b}(\eta) \nu \circ \eta), \nu \circ
\eta \> = 0\, , 
$$
where the derivatives are evaluated at $s=0$.
As $J^{b}(\eta)$ is nowhere zero on $\partial \Omega$,
this is an equation that can be solved for $\partial _{s}q$. This and
(\ref{bd}) complete the desired linearization of (\ref{eq3}). 

We thus obtain the following:

\begin{proposition}
The derivative of $F$ in {\rm (\ref{eqs})} 
at $(\eta , \dot{\eta })$ in the direction of $(u,v)$ is given by
$$
\frac{d}{dt}\left( \begin{array}{c}
                   u \\ \rho v 
                   \end{array}\right)  
=D_{(\eta , \dot{\eta })}F \left( \begin{array}{c}
                                 u \\ v 
                                 \end{array}
                          \right) = 
      \left(
           \begin{array}{c}
           v \\ D_{(\eta , \dot{\eta })}A 
                  \left( 
                       \begin{array}{c}
                       u \\ v 
                       \end{array}
                  \right) 
           \end{array}
      \right) \, , 
$$
where
$$
\begin{array}{rcl}
   D_{(\eta , \dot{\eta })}A\left( \! \!  
                                  \begin{array}{c}
                                       u \\ v
                                  \end{array}
                            \! \! \right) & = &
  A_{ij}^{\alpha \beta}(D\eta) \partial_i \partial_j u^{\beta} +
(\partial_{(\partial _{k} \eta ^{\gamma })}   
 A_{ij}^{\alpha \beta } )(D\eta) 
(\partial _{i}\partial _{j}\eta ^{\beta })(\partial _{k}u^{
\gamma }) +\left[ 
                 \bar{u}\cdot \nabla ,\nabla 
           \right]
        _{\eta }q+ \nabla _{\eta }h \, , \vspace{1mm} \\
\Delta _{\eta }h  & =  &  -{\rm div}_{\eta }(    
  A_{ij}^{\alpha \beta}(D\eta) \partial_i \partial_j u^{\beta} +
(\partial_{(\partial _{k} \eta ^{\gamma })}   
 A_{ij}^{\alpha \beta } )(D\eta) 
(\partial _{i}\partial _{j}\eta ^{\beta })(\partial _{k}u^{
\gamma }) ) - \\ 
& & \left[\bar{u}\cdot \nabla ,{\rm div} \right] _{\eta }{\rm Div}\, W'(D \eta)
-\left[\bar{u}\cdot \nabla ,\Delta  \right] _{\eta }q  + \vspace{1mm} \\ 
    & & 2\rho {\rm trace}(-D\bar{u}(D\eta ^{-1}D\dot{\eta })^{2}+
(D \eta )^{-1}D\bar{v}(D\eta )^{-1}D\dot{\eta }) \; , \vspace{1mm} \\
h\mid _{\partial \Omega }  & = & 
{\displaystyle -\frac{1}{J^{b}(\eta)}\left(  A_{ij}^{\alpha \beta}(D\eta)\partial_j u^{\beta}
N^i \nu^\alpha \circ \eta  +
q \left[ \det{\left(
 \begin{array}{c}
\partial_{T_{1}} u \\ \partial_{T_{2}}\eta 
 \end{array} \right)}+
 \det{\left( \begin{array}{c} 
   \partial_{T_{1}}\eta \\ \partial_{T_{2}} u 
    \end{array} \right)}
   \right] \right)
 } \, , 
\end{array} 
$$
and $q$ solves the boundary value problem {\rm (\ref{eqb})}.
Here $\bar{u}=u\circ \eta ^{-1}$, $\bar{v}=v\circ \eta ^{-1}$,  and
$\{ T_1, T_2\}$ is an orthonormal frame of the boundary.  
\label{eq:jav} 
\end{proposition}

We refer the reader to \cite[Proposition 2.8]{ebsi} or \cite[\S 6]{ebsi2} for
additional details.   

We now derive the linearization of (\ref{eqs}) at the stationary solution
$(\eta(t)(x)=x,0)$. We bear in mind  the information contained in (\ref{eq3}) 
indicating that $W'(D\eta)N$ is normal to the boundary if $\eta$ is a 
solution to (\ref{eq1}), (\ref{eq2}), (\ref{eq3}). We obtain:

\begin{corollary} \label{co4}
The tangent space at the stationary solution $\eta =\BOne$ 
consists of vector fields $u$ such that
\begin{equation} \label{nblg}
\begin{array}{c}
{\rm div}\, u =0 \, ,  \\
(\< \partial _{T}u, W'(\BOne)N \> + \< T,\partial_s(W'(D\eta(s))N)\mid_{s=0} \>
 )\mid_{\partial \Omega} 
= 0\, , 
\end{array} 
\end{equation}
for any vector field $T$ tangent to $\partial \Omega $, and the
the linearization of {\rm (\ref{eqs})} about the pair 
$(\eta (t),\dot{\eta}(t))= (\BOne, 0)$ in the direction of $(u,v)\in \mc{C}_{k,
\alpha}$, $k\geq 1$,  is given by  
$$
\frac{d}{dt}\left( \begin{array}{c}
                   u \\ \rho v 
                   \end{array}\right) = 
      \left(
           \begin{array}{c}
           v \\ A_{ij}^{\cdot \, 
\beta}(\BOne) \partial_i \partial_j u^{\beta}+ \nabla h
           \end{array}
      \right) \, , 
$$
where 
$$
\begin{array}{ccl}
\Delta h  & =  &  0 \\
h\mid _{\partial \Omega }  & = & 
{\displaystyle - A_{ij}^{\alpha \beta}(\BOne )\partial_j u^{\beta}
N^i N^\alpha - \< W'(\BOne)N,N\> \< \partial_N u, N\> \, .}
\end{array} 
$$ 
\end{corollary}

We shall discuss numerical solutions of the associated Cauchy problem for
the system above in the presence of an external force, and for a
body $\Omega \subset \mb{R}^3$ one of whose 
dimensions is one order of magnitude smaller in size than the other two.

Theorem \ref{th2} follows by applying the contraction mapping principle 
\cite{ebsi2}. When the linearization of (\ref{eqs}) is carried out at a point 
$(\eta ,\dot{\eta })$ other than $(\eta(t)x=x,0)$, the resulting 
operator is hyperbolic on the tangent space at $\eta$ of the submanifold 
defined by (\ref{eq1}), a space that depends on $\eta$, and this  
fact makes it difficult to prove the well-posedness of its 
associated Cauchy problem in the usual manner.
 We are forced to enlarge the space where the 
linearization of $F$ is to be defined, and make it independent of the point 
where we carry the linearization. In turn, 
this forces us to modify the equation defined by $F$ to preserve 
ellipticity of the spatial part of its linearization. Theorem \ref{th2} 
follows after a  judicious execution of this technical strategy. The fixed 
point in the contraction mapping that yields a solution to the said modified 
equation is a solution of the unmodified one since, at 
such a point, the modifying term in the equation vanishes \cite{ebsi2}.    

\section{Hooke's law: The stored energy of orthotropic materials}
\label{s3}
Let $n$ be a unit vector at a point $x$ on a elastic body $\Omega$. At $x$, 
the intensity of the contact force per unit area on the plane $\{ y\in 
\Omega: \<y-x  , n\>=0\}$    exerted by the material on $\{ y\in \Omega: \<y-
x  , n\> \geq 0\}$   on that of $\{ y\in \Omega: \<y-x  , n\><0\}$ is 
measured by the stress tensor $T(t)(x,n)$. If the motion is of type $C^1$, 
there exists a $(2,0)$-tensor field $\sigma(t)(x)$ such that $T(t)(x,n)=\< 
\sigma(t)(x),n\>$. In this expression, we are using the metric to think of 
$\sigma(t)(x)$ as a $(1,1)$ tensor.  We call $\sigma$ the stress
tensor of the body. Each of its nine components has the dimension of force 
per unit area, which we think of as pressure functions.
The equilibrium of a body requires the conservation of its momentum. This 
holds if, and only if, the tensor $\sigma $ is symmetric, which we assume 
hereupon. 

The strain tensor is defined measuring the relative change of the metric 
tensor on the body as it changes position in time. If 
$ds_0^2= g^0_{ij}dx^i dx^j$ is the metric at time $t=0$ 
and $ds^2= g_{ij}dy^i dy^j$ is the metric tensor at $y=\eta(t)(x)$, then we 
have $ds^2_0= g^0_{ij}\frac{\partial{x^i}}{\partial y^l} 
\frac{\partial{x^j}}{\partial y^m} dy^i dy^j$, and we obtain that 
$$
ds^2-ds_0^2=\left( g_{ij}-g^0_{ml}\frac{\partial{x^m}}{\partial y^i} \frac{\partial{x^l}}{\partial y^j}
\right) dy^i dy^j=2 e_{ij}dy^i dy^j\, ,
$$ 
defining the strain $2$-tensor as  
$$
e_{ij}=\frac{1}{2}\left( g_{ij}-g^0_{ml}\frac{\partial{x^m}}{\partial y^i} \frac{\partial{x^l}}{
\partial y^j}\right) \, .
$$
It is usually referred to as the Cauchy strain tensor.  
 
In terms of the displacement vector $u$ given by $\eta (t)(x)=x+u(t,x)$, and for the Euclidean metric
in $\mb{R}^3$, we obtain that
\begin{equation} \label{eq13}
e_{ij}=\frac{1}{2}(\nabla u +\nabla u^T - \nabla u^T \nabla u)_{ij} \, .  
\end{equation}
Modulo quadratic errors in $\nabla u$,  
$e$ coincides with the symmetrized vector of covariant derivatives 
of $u$, and we shall often (if not always) equate the two.
Usually this latter notion is called the infinitesimal strain. 
The components of $e$ are dimensionless. 
For instance, if in cylindrical coordinates we write 
the displacement vector as   
$$
u= u_r \partial_r + u_\theta \frac{1}{r}\partial _\theta + u_z \partial_z\, ,
$$ 
then the components of $e$ are 
$$
\begin{array}{lll}
e_{rr}=\partial_r u_r \, , & e_{\theta \theta}= \frac{1}{r}\left( 
\partial_\theta u_\theta +u_r\right) \, , & e_{zz}=\partial_z u_z 
\, , \vspace{1mm} \\
e_{\theta z}=\frac{1}{2}\left( \partial_z u_\theta \! + \!  \frac{1}{r}
\partial_\theta u_z \right)\, ,  & 
e_{zr }=\frac{1}{2}(\partial_r u_z \! + \!  \partial_z u_r)\, , & 
e_{r \theta }=\frac{1}{2}\left(  \frac{1}{r}\partial_\theta u_r \! + 
\partial_r u_\theta \! - \! \frac{u_\theta}{r} \right) \, .  
\end{array}
$$

By computing two covariant derivatives of the tensor $e$, and comparing the
results, we obtain the St. Venant compatibility conditions
\begin{equation}
e_{ij,kl}+e_{kl,ij} -e_{ik,jl}-e_{jl,ik}=0 \, .
\end{equation}

\subsection{Generalized Hooke's law}
Let $\mc{S}^2$ denotes the space of symmetric $2$-tensors.
A body is said to be a generalized Hooke body if there exists a linear operator
$$
W: C^{\infty}(\Omega; \mc{S}^2) \rightarrow C^{\infty}(\Omega; \mc{S}^2)
$$
such that 
\begin{equation}
\sigma = W e \, , 
\end{equation}
and whose stored energy function is given by 
\begin{equation} \label{enc}
W(D \eta)=\frac{1}{2}\< \sigma, e\>=\frac{1}{2}\< We,e\>\, . 
\end{equation}
The tensor $W$, which encodes the material properties of the body,  
is referred to as the tensor of elastic constants, or moduli, of 
the material. 

We use pair of indexes to parametrize symmetric two tensors. In
$3$-space, these tensors have $6$ degrees of freedom. If we express
$W=(W^{ijkl})$ into components, we have that
\begin{equation} \label{dec}
\sigma^{ij}=W^{ijkl} e_{kl} \, ,  
\end{equation}
and
\begin{equation} 
W(D \eta)=\frac{1}{2}\< \sigma, e\>=\frac{1}{2}W^{ijkl} e_{kl} e_{ij}\, , 
\end{equation}
with the symmetries $W^{ijkl}=W^{jikl}=W^{ijlk}$, for $36$ degrees of
freedom in $W$ altogether.  
Since we consider coercive positive definite stored energy functions 
in our work, 
we must have the additional symmetry $W^{ijkl}=W^{klij}$ also, and this 
reduces the degrees of freedom of $W$ down to $21$. 
They are explicitly given by
\begin{equation} \label{ten}
\left( \begin{array}{c}
\sigma_{11} \\
\sigma_{22} \\
\sigma_{33} \\
\sigma_{23} \\
\sigma_{31} \\
\sigma_{12} 
\end{array}
\right)
= 
\left( \begin{array}{cccccc}
W_{11}^{\phantom{1}11} &
W_{11}^{\phantom{1}22} &
W_{11}^{\phantom{1}33} &
W_{11}^{\phantom{1}23} &
W_{11}^{\phantom{1}31} &
W_{11}^{\phantom{1}12} \\
W_{22}^{\phantom{1}11} &
W_{22}^{\phantom{1}22} &
W_{22}^{\phantom{1}33} &
W_{22}^{\phantom{1}23} &
W_{22}^{\phantom{1}31} &
W_{22}^{\phantom{1}12} \\
W_{33}^{\phantom{1}11} &
W_{33}^{\phantom{1}22} &
W_{33}^{\phantom{1}33} &
W_{33}^{\phantom{1}23} &
W_{33}^{\phantom{1}31} &
W_{33}^{\phantom{1}12} \\
W_{23}^{\phantom{1}11} &
W_{23}^{\phantom{1}22} &
W_{23}^{\phantom{1}33} &
W_{23}^{\phantom{1}23} &
W_{23}^{\phantom{1}31} &
W_{23}^{\phantom{1}12} \\
W_{31}^{\phantom{1}11} &
W_{31}^{\phantom{1}22} &
W_{31}^{\phantom{1}33} &
W_{31}^{\phantom{1}23} &
W_{31}^{\phantom{1}31} &
W_{31}^{\phantom{1}12} \\
W_{12}^{\phantom{1}11} &
W_{12}^{\phantom{1}22} &
W_{12}^{\phantom{1}33} &
W_{12}^{\phantom{1}23} &
W_{12}^{\phantom{1}31} &
W_{12}^{\phantom{1}12} 
\end{array}
\right)
\left( \begin{array}{c}
e_{11} \\
e_{22} \\
e_{33} \\
e_{23} \\
e_{31} \\
e_{12} 
\end{array}
\right) \, .
\end{equation}

For a generalized Hooke body $\Omega$,  
if $\eta(s)=\BOne +su$, we have that   
\begin{equation} \label{bs}
\begin{array}{rcl}
W'(\BOne)^{i}_{\mbox{}\hspace{1mm}\alpha} & = & {\displaystyle 
\frac{\partial W}{\partial( \partial_i \eta^\alpha)}\mid_{\eta =\BOne}=
W^{i \hspace{1mm} j}_{\mbox{}\hspace{1mm}\alpha\hspace{1mm}j}} \, , \vspace{1mm}
\\
\partial_s(W'(D\eta(s))N)\mid_{s=0} & = &  
W^{\alpha \hspace{1mm} j}_{\mbox{}\hspace{1.5mm}i \hspace{1.5mm}\beta} 
\partial_j 
u^\beta N^i=\sigma(\nabla u)^{\alpha}_{\mbox{} \hspace{1mm} i}N^i=
\sigma(\nabla u)N \, , \vspace{1mm} \\
A_{ij}^{\alpha \beta}( \BOne ) & = &
W^{i\hspace{1.5mm}j}_{\phantom{i}\alpha \hspace{1.5mm} \beta}\, , 
\end{array}
\end{equation}
and the boundary value of the function $h$ in Corollary \ref{co4} is given by 
\begin{equation}\label{hb}
\begin{array}{rcl}
h\mid_{\partial \Omega} & = &  -N^i N^\alpha W^{\phantom{i}\hspace{1.5mm}j}_
{i\alpha \hspace{1.5mm} \beta} \partial_j u^\beta -  
W^{\phantom{i}\hspace{1.5mm}j}_{i\alpha \hspace{1.5mm} j}N^i N^{\alpha}
\< \partial_N u,N\> \\
& = & -\< \sigma(D u)N, N\>- \< W'(\BOne)N, N\> \< \partial_N u, N\>\, . 
\end{array} 
\end{equation}
Notice that for the isotropic neoHooke materials of earlier, this 
$h\mid_{\partial \Omega}$ is 
the function $-2\< \partial_N u, N\>$, a fact consistent with the case 
analyzed in  \cite{ebsi}.   

Orthotropic materials are generalized Hooke bodies with 
additional properties. They posses three mutually orthogonal planes of 
symmetries at each point, with three corresponding orthogonal axes. 
Their elastic coefficients
are unchanged under a rotation of $180^{\circ}$ about any of these axes, and 
their tensor of elastic constants have nine degrees of 
freedom only. Relative to the preferred axes,  
the tensor of elastic constants (\ref{ten}) reduces to 
\begin{equation}\label{ten2}
W=\left( \begin{array}{cccccc}
W_{11}^{\phantom{1}11} &
W_{11}^{\phantom{1}22} &
W_{11}^{\phantom{1}33} &
0 & 0 & 0  \\
W_{22}^{\phantom{1}11} &
W_{22}^{\phantom{1}22} &
W_{22}^{\phantom{1}33} &
0 & 0 & 0  \\
W_{33}^{\phantom{1}11} &
W_{33}^{\phantom{1}22} &
W_{33}^{\phantom{1}33} &
0 & 0 & 0  \\
0  & 0  & 0 &
W_{23}^{\phantom{1}23} &
0 & 0 \\
0  & 0  & 0 & 0  &
W_{31}^{\phantom{1}31} &
0 \\
0 & 0 & 0 & 0 & 0 &
W_{12}^{\phantom{1}12} 
\end{array}
\right) \, .
\end{equation} 

Wood is normally taken as an example of an orthotropic material.
It has unique and independent mechanical properties along 
three mutually perpendicular directions: The longitudinal axis $z$ that is 
parallel to the fiber grains; the radial axis $r$ that is normal to the growth 
rings; and the tangential axis $\theta$ to the growth rings. 

The elastic constants of an orthotropic Hooke body are described in terms of 
the three moduli of elasticity, the six Poisson ratios, and the three
moduli of rigidity, or shear modulus, determined by the axes of symmetry.
These parameters are the coefficients of
proportionality between the elongation, or compression, experienced by the
material in the direction of a tension, or compression, stress acting
on it, the ratio of the contraction perpendicular to the
direction of the load to the longitudinal strain parallel to it, 
and the ratio of the shear stress to the shear strain when the deflection 
of the material is caused by a shear stress, respectively. 
The moduli of elasticity $e_1, e_2, e_3$, and moduli of rigidity $g_1, g_2,
g_3$, have dimensions of force per unit area, while the 
Poisson ratios $\mu_{12}, \mu_{21}, \mu_{13},\mu_{31}, \mu_{23},
\mu_{32}$ are dimensionless. These twelve parameters satisfy the compatibility
relations
\begin{equation} \label{rel} 
\frac{\mu_{ij}}{e_i}=\frac{\mu_{ji}}{e_j}\, , \quad \text{$i\neq j$, 
$1\leq i,j \leq 3$}\, .
\end{equation}

If we place the vertical axis of a tree trunk along the $z$-direction, 
the axes of symmetry of its wood coincide with the ordered 
cylindrical coordinates $(r,\theta, z)$ referred to earlier. 
Then the elastic behaviour of wood 
is described by the elasticity moduli $e_r, e_\theta, e_z$; the rigidity 
moduli $g_{z r}, g_{z\theta}, g_{r \theta}$; and six Poisson ratios $\mu_{r 
\theta},\mu_{\theta r}, \mu_{r z}, \mu_{z r}, \mu_{\theta z}, \mu_{z \theta}$.
These constants satisfy the three relations (\ref{rel}). 

Values of the constants above for materials relevant to our work 
are given in the tables below. They have been extracted from those given in
\cite{gwk}. The values for $e_z$ displayed were obtained from measurements
carried out on wood with a 12\% of moisture content. We shall use these
in our simulations here. Accurate values for wood specifically used in the 
construction of actual plates should lead to results more faithful to the 
reality under analysis.
\medskip

\begin{center}
\begin{tabular}{|l|c|c|c|c|c|c|} \hline \hline
  & $e_z$ in MPa & $e_{\theta}/e_z$ & $e_{r}/e_z$ & $g_{z r}/e_z$ &  
 $g_{z\theta}/e_z$ & $g_{r \theta}/e_z$  \\
\hline \hline
Maple, sugar & 13,860 & 0.065 & 0.132 & 0.111 & 0.063  & ---   \\
\hline
Maple, red   & 12,430 & 0.067 & 0.140 & 0.133  & 0.074  & --- \\
\hline
Spruce, Sitka & 10,890 &  0.043 & 0.078 & 0.064  & 0.061  & 0.003  \\
\hline
Spruce, Engelmann & 9,790 & 0.059  & 0.128  & 0.124  & 0.120  & 0.010  \\
\hline \hline 
\end{tabular}
\medskip 

\centerline{Table 1. Ratios of elasticity to rigidity moduli for spruce and 
maple.}
\end{center}

\mbox{}
\medskip
\begin{center}
\begin{tabular}{|l|c|c|c|c|c|c|} \hline \hline
  & $\mu_{z r}$ & $\mu_{z\theta}$ & $\mu_{r \theta}$ & $\mu_{\theta r}$ & $\mu_{r z}$ & 
$\mu_{\theta z}$  \\
\hline \hline
Maple, sugar & 0.424 & 0.476 & 0.774 & 0.349  & 0.065 & 0.037  \\
\hline
Maple, red   & 0.434 & 0.509 & 0.762  & 0.354  & 0.063 & 0.044 \\
\hline
Spruce, Sitka & 0.372 & 0.467 & 0.435  & 0.245 & 0.040 & 0.025  \\
\hline
Spruce, Engelmann & 0.422  & 0.462 & 0.530 & 0.255 & 0.083 & 0.058  \\
\hline \hline 
\end{tabular}
\medskip

\centerline{Table 2. Poisson ratios for spruce and maple.}
\end{center}
\medskip

The inverse of the tensorial relation (\ref{dec}) yields
\begin{equation} \label{decin}
e^{ij}= U^{ijkl}\sigma_{kl} \, , 
\end{equation}
where $U=W^{-1}$. For orthotropic materials, and in the ordered cylindrical 
coordinates $(r,\theta, z)$ above, we must have that 
\begin{equation} \label{decI}
U=W^{-1}= \left( \begin{array}{cccccc}
{\displaystyle \frac{1}{e_r} } &
{\displaystyle -\frac{\mu_{\theta  r}}{e_\theta}}  &
{\displaystyle -\frac{\mu_{z r}}{e_z}}  &
0 & 0 & 0  \vspace{1mm} \\
{\displaystyle -\frac{\mu_{r \theta }}{e_r}} &
{\displaystyle \frac{1}{e_\theta}}  &
{\displaystyle -\frac{\mu_{z \theta}}{e_z}}  &
0 & 0 & 0  \vspace{1mm} \\
{\displaystyle -\frac{\mu_{r z}}{e_r}} &
{\displaystyle -\frac{\mu_{\theta z}}{e_\theta }} &
{\displaystyle \frac{1}{e_z}} &
0 & 0 & 0  \vspace{1mm} \\
0  & 0  & 0 &
{\displaystyle \frac{1}{g_{\theta z}} } &
0 & 0 \vspace{1mm} \\
0  & 0  & 0 & 0  &
{\displaystyle \frac{1}{g_{rz}}} &
0 \vspace{1mm} \\
0 & 0 & 0 & 0 & 0 &
{\displaystyle \frac{1}{g_{r\theta} }} 
\end{array}
\right) 
\, .
\end{equation}
For suppose the body is subjected to simple tension or compression with 
tension stress in the radial direction only. The only nonzero component of 
the tensor $\sigma$ is $\sigma_{rr}$, and by (\ref{decin}) and (the first 
column of) (\ref{decI}), we must have that
$$
e_{rr} =\frac{1}{e_r}\sigma_{rr}\, , \quad 
e_{\theta \theta } =-\frac{\mu_{\theta r}}{e_\theta}\sigma_{rr}\, , \quad 
e_{zz} =-\frac{\mu_{z r}}{e_z}\sigma_{rr}\, , 
$$ 
respectively. It follows that $e_r$ is the quotient of $\sigma_{rr}$ and 
$e_{rr}$, which identifies this constant as the radial modulo of elasticity.
By (\ref{rel}), we then obtain that 
$e_{\theta \theta}=-\mu_{\theta r}e_{rr}$ and
$e_{zz}=-\mu_{z r}e_{rr}$, and so 
$\mu_{\theta r}$ and $\mu_{zr}$ are the 
corresponding Poisson ratios.  
The form of the remaining coefficients in the second and third columns 
of $U$ are obtained similarly by subjecting the body to simple tensions or 
compressions in the angular or longitudinal directions instead.
Finally, if we consider the case of pure shear with the only nonzero 
components of $\sigma$ being $\sigma_{r \theta}=\sigma_{\theta r}$, proceeding 
as above, we conclude that
$$
e_{r\theta }=\frac{1}{g_{r\theta}} \sigma_{r\theta } \, ,
$$
and so $g_{r\theta}=\sigma_{r\theta}/e_{r\theta }$ is identified as the
shear modulo $g_{r\theta}$.  
Similar considerations allow for the interpretation of the remaining 
diagonal shear modulo terms in (\ref{decI}).

\begin{remark} \label{r6}
The matrix of change of coordinates from the orthonormal
basis $\{ \partial_r, \frac{1}{r}\partial_\theta, \partial_z \}$ to the 
orthonormal basis $\{ \partial_x, \partial_y, \partial_z\}$ is given by
$$
\left( \begin{array}{crc}
\cos{\theta}  & -\sin{\theta} & 0 \\
\sin{\theta}  & \cos{\theta}  & 0 \\
 0 & 0\phantom{bb} & 1 
\end{array} 
\right)\, .
$$
On wedges $\{ (r,\theta, z): \; r\leq r_0 \, ,  |\theta| \leq
\varepsilon/2 \}$ of small total angular variation $\varepsilon$, this 
matrix is very close to the identity. On these  
regions, and for this reason, the components of the tensor $W$
of an orthotropic material in cylindrical and Cartesian coordinates will 
be considered to be the same.          
\end{remark}

\section{Numerical method} \label{s4}
We consider the linear equation of Corollary \ref{co4} in the 
presence of an external sinusoidal force $F_{\omega}=F_{\omega}(t,x)$, and 
with trivial initial data:
\begin{equation}\label{vp}
\begin{array}{lcl}
{\displaystyle \rho \frac{d^{2}u}{dt^{2}}} & = & A_{ij}^{\cdot \, \beta}(\BOne)
\partial_i \partial_j u^\beta + \nabla h +F_\omega \, , \\
u\mid_{t=0} & = & 0 \, , \\
\partial_t u\mid_{t=0} & = & 0 \, .
\end{array}
\end{equation}
Here, $h$ is the harmonic extension to $\Omega$ of the function
(\ref{hb}) over $\partial \Omega$:
$$
h=\mc{H}( -\< \sigma (D u)N,N\>- \< W'(\BOne)N, N\>\< \partial_N u, N\>)
\, .
$$
We seek a numerical solution to this system that satisfies the conditions in 
(\ref{nblg}). We describe here the algorithm we use for that purpose, and 
how we apply it to determine the vibration modes of the body in motion.  

\subsection{Divergence free vector fields associated to smooth triangulations} 
For simplicity of the exposition, we work first in the smooth category.
Let $K$ be a finite smooth oriented triangulation of an oriented connected 
manifold 
with boundary $M$ that has been provided with a fixed Riemannian metric $g$. 
In our work, $g$ is always the metric induced by the Euclidean metric in 
$\mb{R}^3$, where our body $\Omega$ is embedded. 
We identify the polytope of $K$ with $M$, and fix some ordering of the 
vertices of $K$. We denote by $C^q(K)$ the space of simplicial {\it oriented} 
$q$-cochains, and by $L^2(M;\Lambda^q T^*M)$ the space of $L^2$ $q$-forms on 
$M$. If $P \in C^{\infty}(M; E) \rightarrow C^{\infty}(M;F)$ is a linear 
operator over $M$ mapping sections of a bundle $E$ to sections of a bundle $F$,
we denote by $L^2_P(M; E)$ the subspace of $L^2(M;E)$ sections that are mapped
by $P$ into $L^2(M;F)$, provided with the graph norm.     
  
Every cochain $c$ can be written uniquely as $c=\sum c_\sigma \sigma^*$, 
where the sum is taken over all simplices $\sigma=[p_0, \ldots, p_q]$ of
$K$ whose vertices form an increasing sequence with respect to the ordering, 
and where the $c_{\sigma}$s are real coefficients. In this expression, 
$\sigma^*$ is 
the characteristic function ordered cochain determined by the simplex $\sigma$. 

For a vertex $p$, we let $x_p$ be the $p$-th barycentric coordinate in $K$, a 
continuous piecewise linear function with support in $\overline{\rm St}\, p$, 
the closure of the star of $p$ in $K$. If a point $x$ belongs to a simplex
spanned by a set of vertices that includes $p$, $x_p(x)$ is the weight of
$p$ in the convex combination that expresses $x$ as a weigthed sum of the
said vertices. The collection of functions $\cup_p \{ x_p \}$ forms a 
continuous partition of unity of the underlying space, or polytope, of $K$,
subordinated to the open cover $\{ {\rm St}\, p\}$.
       
The barycentric subdivision $K'$ of a (not necessarily oriented) triangulation
$K$ is a simplicial complex that is naturally oriented. 
For if $\sigma$ is a simplex in $K$, we denote by 
$\hat{\sigma}$ its barycenter in $K'$.
For any two simplices $\sigma_i, \sigma_j$ of $K$, the expression 
$\sigma_i \succ \sigma_j$ 
means that $\sigma_j$ is a proper face of $\sigma_i$. Then the simplices of 
$K'$ are of the form
$$
[\hat{\sigma}_{i_1},\hat{\sigma}_{i_2}, \ldots, \hat{\sigma}_{i_k}]
$$
where $\sigma_{i_1}\succ \sigma_{i_2} \succ \cdots \succ \sigma_{i_k}$. 
The vertices of $K'$ are ordered by decreasing dimension of the 
simplices of the triangulation of which they are the barycenters. This 
ordering induces a linear ordering of the vertices of each simplex of $K$
\cite{mu}.  

In what follows, we shall consider 
the barycentric subdivision triangulation $K'$ of a given smooth 
triangulation $K$ of a Riemannian $3$-manifold $(M,g)$. 
We shall denote by $K^{(j)}$ the $j$th skeleton of $K'$,  
and by $\mc{V}'$, $\mc{E}'$, $\mc{F}'$, and $\mc{T}'$, the set of vertices, 
edges, faces, and tetrahedrons   
of $K'$, respectively. The {\it interior edges}, those edges in $K^{(1)}$ 
with at least one of its vertices lying in the interior of $M$,
are denoted by $\mc{E}'_{\circ}$, while the {\it boundary edges}, those 
whose vertices are contained on the boundary of $M$, are denoted by 
$\mc{E}'_{\partial }$. The {\it interior faces}, those where at least one of
its three vertices is contained in the interior of $M$, shall be denoted by
$\mc{F}'_\circ$, while the {\it boundary faces}, those whose three vertices are 
on the boundary of $M$, shall be denoted by $\mc{F}'_\partial$. 
When necessary, the analogous concepts for the triangulation $K$ itself will 
be denoted similarly but without the $'$s. Notice that $K^{(1)}$ is the 
oriented graph $G(\mc{V}', \mc{E}')$. 
The cardinality of a set $S$ will be denoted by $|S|$.  

We can attach to the oriented triangulation $K'$ a natural 
set of piecewise linear forms, the Whitney forms \cite{whi}, and their
metric dual vector fields. We derive discrete approximations to 
the solutions of (\ref{vp}) by using the Whitney vector fields associated to 
the faces of $K'$, but it is the capturing of the algebraic topology of $M$ by
the entire collection of Whitney forms that is key in explaining why 
these approximate solutions are the best to be considered in our setting.
For general definitions and properties of the Whitney forms, we refer the 
reader to \cite{do}.  
  
Given any oriented edge $e=[p_0,p_1]$ in $\mc{E}'$, we consider the 
piecewise continuous $1$-form
$w_e= x_{p_0} d x_{p_1}- x_{p_1}dx_{p_0}$, and its metric dual vector field
\begin{equation} \label{f1e}
W_e= x_{p_0} \nabla^g x_{p_1}- x_{p_1}\nabla^g x_{p_0} \, ,
\end{equation}
as elements of the space of $L^2$ forms or vector fields with integrable 
square norm. The vector field $W_e$ is compactly supported in 
$\overline{\rm St}\, p_0 \cap \overline{\rm St}\, p_1=\overline{{\rm St}}\, e$.
It is divergence-free in the interior of any symplex of maximal dimension in 
its support, but has singular distributional normal derivative on codimension 
$1$ simplices in the boundary of these, although its tangential derivatives
are smooth.  If $e'$ is any oriented edge of $K'$, we have that
\begin{equation} \label{ee}
e^*(e')=\< w_e, e'\> =\int_{e'} W_e = \left\{ \begin{array}{ll}
                      0 & \text{if $e'\neq e$}\, , \\
                      1 & \text{otherwise}\, .
                      \end{array}
              \right.    
\end{equation}

The span of $\{ W_e\}_{e\in \mc{E}'}$ contains many gradients 
vector fields.
 Indeed,  
let $b_{pe}$ be the incidence number of the vertex $p$ and edge $e$ in the
oriented graph $G(\mc{V}',\mc{E}')=K^{(1)}$. This number is $1$ if
$e=[q,p]$, $-1$ if $e=[p,q]$, and zero otherwise. Therefore,
$$
\sum_{e\in \mc{E}'} b_{pe}W_e = \sum_{q\in \mc{V}'} (x_q \nabla^g x_p -
x_p\nabla^g x_q)=(\sum_{q\in \mc{V}'} x_{q})\nabla^g x_p -
x_p \nabla^g ( \sum_{q\in \mc{V}'} x_q )   =\nabla^g x_p \, . 
$$

\begin{lemma}
The set $\{ W_e \}_{\text{$e=[p_0,p_1]\in \mc{E}'$}}$ is an orthonormal 
linearly independent set of $L^2$ vector fields. Its span contains the 
span of the set of gradients $\{ \nabla^g x_p \}_{p \in \mc{V}'}$.   
\end{lemma}

The space ${\rm span} \{ W_e \}_{\text{$e=[p_0,p_1]\in \mc{E}'$}}$ is thus
suitable for discretizing gradient vector fields.   
In order to discretize divergence-free fields instead,  
we should look for a space of vector fields that contains many curls.
We proceed as follows. 

We first recall that if $d\mu $ stands for the 
standard orientation form, the 
Hodge start operator $* : \Lambda^p T^*M \rightarrow \lambda^{3-p}T^*M$ is 
defined uniquely by the identity
$$
\alpha \wedge \beta = \< *\alpha, \beta \> d\mu \, , 
$$    
and produces an isomorphism of $C^\infty(M)$ algebras. We make use of the
correspondence between $1$ and $2$ forms 
that $*$ defines.

Let $f=[p_0,p_1,p_2] \in \mc{F}'$ be a face in $\mc{F}'$. The Hodge $*$ of
the $L^2$ form $w_f= 2(x_{p_0}dx_{p_1} \wedge d x_{p_2}+x_{p_1}dx_{p_2} 
\wedge d x_{p_0} +x_{p_2}d x_{p_0} \wedge d x_{p_1})$
is an $L^2$ one form with metric dual  
\begin{equation} \label{f2f}
W_f= 2(x_{p_0}\nabla^g x_{p_1} \times \nabla^g x_{p_2}+x_{p_1}\nabla^g x_{p_2} 
\times \nabla^g x_{p_0} +x_{p_2}\nabla^g x_{p_0} \times \nabla^g x_{p_1}) \, .
\end{equation}
Here, $\nabla^g x_{p_i} \times \nabla^g x_{p_j}$ stands for the cross product
of the piecewise constant 
$L^2$-gradients $\nabla^g x_{p_i}$ and $\nabla^g x_{p_j}$, 
respectively. The flux of $W_f$ through a face $f'$ in $\mc{F}'$ is 
well-defined, and given by  
\begin{equation} \label{ff}
f^*(f')=\int_{f'} W_f \cdot n_{f'}  = \left\{ \begin{array}{ll}
                      0 & \text{if $f'\neq f$}\, , \\
                      1 & \text{otherwise}\, .
                      \end{array}
              \right.
\end{equation}

Let $e=[p_0,p_1]$ be an edge in $\mc{F}'$, and let  
$W_e$ be the vector field (\ref{f1e}) associated to it, with 
 support on $\overline{{\rm St}}\, e$. Since the tangential derivatives of 
$W_e$ are smooth, its curl is a well-defined piecewise constant vector 
field in $L^2$, and if $b_{ef}$ is now the incidence number of the edge 
$e$ on the face $f$, we have that 
\begin{equation} \label{vf}
U_e:= {\rm curl}\, W_e= 2 \nabla^g x_{p_0}
\times \nabla^g x_{p_1} = \sum_{f\in \mc{F}'} b_{ef}W_f \, .
\end{equation}
Hence, the divergence of $U_e$ in $L^2$ is also well-defined, and equals zero. 

\begin{lemma}
The set $\{ W_f \}_{\text{$f\in \mc{F}'$}}$ is a 
linearly independent set of continuous $L^2$ vector fields whose
span contains the image under {\rm curl} of the span of 
$\{ W_e \}_{e \in \mc{E}'}$.
\end{lemma}

The cohomology of $M$ arises from the homology of either one of the row
complexes in the commutative diagram
$$
\minCDarrowwidth10pt\begin{CD}
0 @>>> C^{\infty}(M) @>{\rm grad}>> C^\infty(M,TM) @>{\rm curl}>> 
C^\infty(M,TM) @>{\rm div}>>  C^\infty(M)  @>>> 0 \\
@.   @AiAA                      @A\sharp AA     @A\sharp *AA      @A*AA    @. \\
0 @>>> C^{\infty}(M) @>d>> C^\infty(M,T^*M) @>d>> C^\infty(M,\Lambda^2 T^*M) @>d>>
C^\infty(M,\Lambda^3 T^*M) @>>> 0  
\end{CD}
$$
which are the same. Given a triangulation of $M$, the vertical arrow maps 
in this diagram identify the $0$ cochain $\sum_p  c_p p^*$ with the function 
$\sum c_p x_p$, the $1$ cochain $\sum_e  c_e e^*$ with the vector field 
$\sum c_e W_e$, the $2$ cochain $\sum_f  c_f f^*$ with the vector field 
$\sum c_f W_f$, and the $3$ cochain $\sum_t c_t t^*$ with the function
$\sum c_t x_{t_b}$, $t_b$ the barycenter of $t$. Thus,
the $2$ cocycles may be viewed as the vector fields $\sum c_f W_f$ such 
that ${\rm div} \sum c_f W_f =0$, while the $2$ boundaries are the set
of elements of the form ${\rm curl}\sum_e c_e W_e$, themselves cocycles 
in their own right. 
The second cohomology of $M$ measures how much larger the $2$ cocycles are
from the $2$ boundaries.

In the cases of interest to us, the polytope of $K$, which we identify with 
$M$, is contractible to a point,
 or to the wedge of two circles. Then the kernel of the divergence operator
on ${\rm span}\, \{ W_f \}_{f\in \mc{F}'}$ coincides with 
the image of ${\rm span} \{ W_e \}_{e \in \mc{E}'}$ under the curl operator.
In general, this is true only modulo a finite dimensional space. 
If $b_{ft}$ is the incidence number of the face $f$ in the tetrahedron $t$, 
we have that 
$$
{\rm div} \, (\sum_{f\in \mc{F}'}c_f W_f) =\sum_{t\in \mc{T}'}
(\sum_{f \in \mc{F}'} c_f b_{ft}) t^{*}\, , 
$$ 
and so, $\sum_{f\in \mc{F}'}c_f W_f$ is divergence-free if, and only if, the 
weighted sum $\sum_{f\in \mc{F}'}c_f b_{ft}$ over the four faces of each
tetrahedron $t$ in $K'$ is identically zero.

Any solution to (\ref{vp}) that satisfies the first of the conditions in
(\ref{nblg}) yields a continuous path in the second cohomology group of
the body, which, therefore, must be constant. A discretization of
the equation should be carried out in a space that preserves this cohomology
element that the solution represents. When doing so by using the Whitney
vector fields $\{ W_f \}_{f\in \mc{F}'}$, the PL nature of the
vector fields we use imposes some obstacles that in practice
we overcome tacitly by exploiting the well-posedness of the equation.

We denote the spaces spanned by the $W_f$s in (\ref{f2f}) and by the $U_e$s in
(\ref{vf}) as
$$
\begin{array}{rcl}
{\rm Div}(K) & = & {\rm span }\{ W_f\}_{f\in \mc{F}'} \, , \\
{\rm Div}_0(K) & = & {\rm span}\{ U_e \}_{\text{$e\in \mc{E}'$}} \, .
\end{array}
$$ 
Their dimensions are  ${\rm dim}\, {\rm Div}(K)=| \mc{F}'|$ and
${\rm dim} \, {\rm Div}_0(K)= | \mc{E}'|$, respectively. 

A weak formulation of problem (\ref{vp}) may be carried out naturally in the
space $L_{\rm div}^2(M;TM)$, or better yet, in its subspace
$$
L_{{\rm div}_0}^2(M;TM) = \{ X \in L^2_{\rm div}(M;TM): \;   
{\rm div}X =0\} \, . 
$$
Accordingly, we could attempt to discretize the said problem 
in either ${\rm Div}(K)$ or ${\rm Div}_0(K)$. But the former is a better
choice. For the vector fields $U_e$s are piecewise constant, and so any of 
their weak derivatives would be trivial on the full measure open subset of 
their support where they are smooth. We thus choose to discretize the weak 
formulation of (\ref{vp}) in (a subspace of) ${\rm Div}(K)$ below.
 
The choice we make raises an issue that is of interest to discuss. 
The nonlinear problem of departure is well-posed in $C^{2,\alpha}$ spaces
for bodies that are of type at least $C^{1, \alpha}$, and so the 
linearized equations are well-posed in $C^{1, \alpha}$. Our discretizing 
space consists of PL vector fields that are quite suitable for maintaining the 
cohomological condition of being divergence free, and perhaps capture
$C^0$ properties of the solution also. But these vector fields are not 
even $C^1$, and so the question arises as if these numerical solutions we 
construct are consistent with the sought after true $C^{1,\alpha}$ solution.
Though a problem we cannot resolve at present by exploring the behaviour 
under mesh refinement due to the high computational complexity, 
the ensuing numerical scheme for finding solutions to the initial value
problem for (\ref{vp}) when discretizing over ${\rm Div}(K)$ should be quite
accurate. If initially the conditions in (\ref{nblg}) hold, the
well-posedness of the equation should force these conditions to hold
also at later times within a small margin of error. And even when extending
the equations, and algorithm, to cover the case of elastic smooth bodies 
with corners, to which all of the spaces above associated to the 
triangulation $K$, as well as the $L^2$ spaces that we considered, have 
natural extensions also, the numerical solutions obtained should 
still reflect some what can be saved of the consistency 
issue, and significant continuity properties of the motion
of these singular elastic bodies, while they vibrate in the elastic regime.

We continue our work often relaxing the $C^{1,\alpha}$ assumption on $M$ to
that of being a smooth manifold with corners. All of the spaces above 
associated  to the triangulation $K$, as well as the $L^2$ spaces that we 
considered, have natural extensions to that context if some of the differential 
operators involve in their definition are interpreted weakly.

\subsection{The discretizing space}
We use the decomposition $\mc{F}'=\mc{F}'_{\circ} + \mc{F}'_{\partial}$ to 
split the representation of an element $U\in {\rm Div}_0(K)\subset {\rm Div}(K)$
into blocks accordingly, 
\begin{equation} \label{asb}
U= \sum_{f\in \mc{F}'} c_f W_f= \sum_{f\in \mc{F}'_{\circ}} c_f W_f+ 
\sum_{f\in \mc{F}'_{\partial }}c_f W_f \, . 
\end{equation}
In the spirit of Einstein summation convention, we express this splitting 
succinctly as $U= c_{f_\circ}W_{f_\circ}+ c_{f_\partial}W_{f_\partial}$.  

In order for $U$ to satisfy a discretized version of the boundary condition 
in (\ref{nblg}) also, we impose a set of $2| \mc{F}_\partial|$-homogeneous 
linear 
equations over the $|\mc{F}'_\partial|$-boundary coefficients of its elements,
 as follows: For each face $f_\partial$ in $\mc{F}_\partial$, let 
$\{T_1, T_2, N\}$ be an oriented orthonormal frame, with $N$ the exterior 
normal to the said face in $K$. Since the boundary conditions are local, 
over this face they involve only the coefficients of $U$ that
are associated to the $6$ faces $f^j_{f_\partial}$, $j=1, \ldots, 6$, in its 
barycentric subdivision.  
We require that
\begin{equation} \label{e27}
\begin{array}{rcl}
\sum_{j=1}^{6} (\< \partial _{T_1}W_{f^j_{f_\partial}}, W'(\BOne)N\>+\< T_1,
\sigma(\nabla W_{f^j_{f_\partial}})N\>) c_{f^j_{f_\partial }} & =  & 0 \, , \\  
\sum_{j=1}^{6} (\< \partial _{T_2}W_{f^j_{f_\partial}}, W'(\BOne)N \> + 
\< T_2, \sigma(\nabla W_{f^j_{f_\partial}})N \>) c_{e^j_{f_\partial }} & = & 
0  \, ,  
\end{array}
\end{equation}
where, by (\ref{bs}), we have that 
$$
\begin{array}{rcl}
(\sigma(\nabla W_{f^j_{f_\partial}}) N)^\alpha & = & 
W^{\alpha \hspace{1mm} k}_{\mbox{}\hspace{1.5mm}i \hspace{1.5mm}\beta}
\, \partial_k W_{f^j_{f_\partial}}^\beta N^i \, , \\
(W'(\BOne)N)^{\alpha}  & = & {\displaystyle
W^{\alpha \hspace{0.5mm} j}_{\mbox{}\hspace{1.4mm}i \hspace{1mm}j}N^i} \, . 
\end{array}
$$
The pairings in these equations are in the sense of $L^2$ of the boundary.

We obtain in this manner an undetermined system of 
$2|\mc{F}_\partial |$ equations in the $6| \mc{F}_\partial|$ 
boundary coordinates of $U$. Since the coupling occurs only among the six 
coefficients associated to the barycentric subdivision faces of a given 
face $f_{\partial }\in \mc{F}_\partial$, the row reduction of the 
associated matrix over these coefficients is a $2\times 6$ row echelon 
form matrix whose rank is either two or one, generically the former.
By a suitable reordering of the basis elements, 
we may write the nonzero rows of the row reduced matrix of the entire system as 
$(\BOne \; -C)$, where $C$ is a block whose number of rows and columns 
are bounded above by $2| \mc{F}_\partial |$ and $5| \mc{F}_\partial|$, 
respectively. We decompose the boundary faces in $\mc{F}'$ accordingly, 
$\mc{F}'_\partial = \mc{F}'_{\partial^I}
+\mc{F}'_{\partial^B}$, so that a vector field
$c_{f_\circ}W_{f_\circ}+c_{f_\partial}W_{f_\partial}=
c_{f_\circ}W_{f_\circ}+
c_{f_{\partial^I}}W_{f_{\partial^I}}+c_{f_{\partial^B}}W_{f_{\partial^B}}$ in
${\rm Div}(K)$ satisfies (\ref{e27}) if, and only if, 
$c_{f_{\partial^I}}=Cc_{f_{\partial^B}}$.
We define ${\rm Div}^{\rm b}(K)$ as such a subspace of ${\rm Div}(K)$: 
\begin{equation} \label{e30bis} 
{\rm Div}^{\rm b}(K)  = \{ c_{f_{\circ}} W_{f_{\circ}}+ 
c_{f_{\partial^I}} W_{f_{\partial^I}}+ c_{f_{\partial^B}} 
W_{f_{\partial^B}}: \; c_{f_{\partial^I}}=C c_{f_{\partial^B}}\} \, .
\end{equation}
If the row reduced matrix of the system of boundary conditions
(\ref{e27}) were not to have any null row, we would have that  
${\rm dim}\, {\rm Div}^{\rm b}(K)=|\mc{F}'|-2| \mc{F}_\partial |=
|\mc{F}'_{\circ}|+4| \mc{F}_{\partial}|$. Otherwise we have that 
${\rm dim}\, {\rm Div}^{\rm b}(K)=|\mc{F}'|-\left( 2(| \mc{F}_\partial |-r_n)+
r_n\right)= |\mc{F}'_{\circ}|+4| \mc{F}_{\partial}|+r_n$, where $r_n$ is the
number of null rows. If $C=(c_{f_{\partial^I},f_{\partial^B}})$, the set
$\{ W_{f_\circ}, c_{f_{\partial^I},f_{\partial^B}}W_{f_{\partial^I}} + 
W_{f_{\partial^B}}\}_{f_\circ \in \mc{F}'_\circ, f_{\partial^B}\in \mc{F}'_{
\partial^B}}$ is a basis for ${\rm Div}^{\rm b}(K)$.

\subsection{The discretized equation}
We discretize a weak solution $U$ of (\ref{vp}) over ${\rm Div}(K)$   
in terms of the vector fields $\{ W_f \}_{f\in \mc{F}'}$. By
the decomposition $\mc{F}'=\mc{F}'_{\circ} + \mc{F}'_{\partial}$, we split $U$
into blocks, 
$$ 
U\sim c_{f_\circ} W_{f_\circ }+ c_{f_{\partial}} W_{f_{\partial}} \, . 
$$ 
and find the components of the vector $\left( \begin{array}{c}
        c_{f_\circ} \\ c_{f_{\partial}} 
        \end{array} 
\right)$ 
by solving the second order differential equation 
that results from the weak formulation of the equation.

 For convenience, we 
use an {\it orthonormal} frame to write the components $A_{ij}^{\alpha \beta}$ 
of the tensor $A$, and assume that this tensor is {\it covariantly constant} 
over the support of each of the basis vectors in ${\rm Div}(K)$. 
If over dots stand 
for time derivatives, we have that
\begin{equation} \label{e24}
\rho I_{\circ, \partial}
 \left( \begin{array}{c}
\ddot{c}_{f'_\circ} \\ \ddot{c}_{f'_{\partial}} 
\end{array} 
\right)= K_{\circ , \partial }
 \left( \begin{array}{c}
c_{f'_\circ} \\ c_{f'_{\partial}} 
\end{array} 
\right) +
\left( \begin{array}{c}
\< F_\omega, W_{f_\circ}\>  \\ \< F_\omega, W_{f_{\partial}} \> 
\end{array} 
 \right) \, ,
\end{equation}
where the matrices $I_{\circ , \partial }$ and $K_{\circ , \partial }$
in this system are 
$$
I_{\circ , \partial}
= \left(  \begin{array}{cc}
I_{f_\circ,f'_\circ} &  I_{f_\circ,f'_{\partial}} \\
I_{f_{\partial},f'_\circ} & I_{f_{\partial},f'_{\partial} }
\end{array}
\right)=
\left( \begin{array}{cc}
\<W_{f_\circ},W_{f'_\circ}\> &  \< W_{f_\circ},W_{f'_{\partial}}\> \\
\< W_{f_{\partial}},W_{f'_\circ}\> & \<W_{f_{\partial}},W_{f'_{\partial} }\>
\end{array}
\right) \, ,
$$
and 
$$
K_{\circ, \partial }= \left(  \begin{array}{cc}
K_{f_\circ,f'_\circ} &  K_{f_\circ,f'_{\partial}} \\
K_{f_{\partial},f'_\circ} & K_{f_{\partial},f'_{\partial} } 
\end{array} 
 \right)=  \left(  \begin{array}{rcl}
-\< \partial_i W^{\alpha}_{f_\circ}, A^{\alpha \beta}_{ij}\partial_j 
W^\beta_{f'_\circ}\> & -\< \partial_i W^{\alpha}_{f_\circ},
A^{\alpha \beta}_{ij} \partial_j W^\beta_{f'_{\partial}}\>  \\
 -\< \partial_i W^\alpha_{f_{\partial}},
A^{\alpha \beta}_{ij}\partial_j W^\beta_{f'_\circ}\> & 
 -\< \partial_i W^\alpha_{f_{\partial}},
A^{\alpha \beta}_{ij} \partial_j W^\beta_{f'_{\partial} }\>  +
B(W_{f_{\partial}}, W_{f_{\partial}'})
\end{array}
\right) \, ,
$$
respectively, the boundary term  
$B(W_{f_{\partial}}, W_{f_{\partial}'})$ 
given by
$$
B(W_{f_{\partial}}, W_{f_{\partial}'}) = 
\< \sigma(D W_{f_{\partial}'})N,W_{f_{\partial}}\> -
(\< \sigma(D W_{f_{\partial}'})N,N\> +\< W'(\BOne) N,N\>\<\partial_N 
W_{f_{\partial}'}, N\>)\< W_{e_{\partial}},N\> \, .    
$$

The matrices in (\ref{e24}) are sparse. Their entries are zero if
${\rm St}\, f \cap {\rm St}\, f'= \emptyset$. In fact, 
all of $\< W_{f_{\partial}}, W_{f'_{\partial} }\>$,
$-\< \partial_i W^\alpha_{f_{\partial}},
A^{\alpha \beta}_{ij} \partial_j W^\beta_{f'_{\partial} }\>$, as well as
the boundary term $B(W_{f_{\partial}}, W_{f_{\partial}'})$ are diagonal, 
that is to say, identically zero if $f_{\partial}\neq f'_\partial$. This
latter fact, and the symmetries of the tensor $A^{\alpha \beta}_{ij}$, 
makes of (\ref{e24}) a symmetric system.

The numerical approximation to our solution of (\ref{vp}) is the vector field
\begin{equation} \label{e30}
U= c_{f_\circ} W_{f_\circ}+ c_{f_{\partial^I}}W_{f_{\partial^I}}+
c_{f_{\partial^B}} W_{f_{\partial^B}} \, ,
\end{equation}
obtained in this manner, with the coefficients 
$\left( \begin{array}{c}
        c_{f_\circ} \\ c_{f_{\partial}}
        \end{array}
\right)$
given by the solutions of (\ref{e24}).

The eigenvalues of $-(\rho I_{\circ, \partial})^{-1} K_{\circ,
\partial}$ and the corresponding frequencies they induce approximate the 
natural vibration frequencies of the body. Positive eigenvalues lead to 
vibrations that decay exponentially fast in time. Negative eigenvalues lead
to undamped vibration modes. By resonance, any one of these 
induces  oscillatory motions within the body when this is subjected to an 
external pressure wave of frequency close to the frequency of the wave
mode determined by the eigenvalue.     

The eigenvectors of $-(\rho I_{\circ, \partial})^{-1} K_{\circ, \partial}$  
do not necessarily satisfy the discrete version of 
either of the two conditions in (\ref{nblg}) (although the well-posedness of 
the initial value problem for (\ref{vp}) should make the discrepancy 
between these and the conditions the eigenvectors satisfy relatively small).
It is for this reason that we take the waves induced by the eigenvectors of
this matrix as {\it coarse} approximations to the vibration patterns of the 
body, and call {\it coarse resonance waves} the waves produced by resonance 
for frequencies close to the frequencies of these coarse approximations.
In our simulations, we shall describe some of these coarse resonance 
waves by depicting their nodal points over the portion of the
boundary opposite to that where the external wave hits it.

We derive {\it fine} approximations by incorporating the boundary conditions 
of (\ref{nblg}) into our solution
of (\ref{e24}), discretizing the solution $U$ over 
the space ${\rm Div}^{\rm b}(K)\subset {\rm Div}(K)$ instead.
We use the natural basis for this space that the splitting 
$\mc{F}'_{\partial}= \mc{F}'_{\partial^I}+ \mc{F}'_{\partial^B}$
leads to, and express the discretized solution into blocks accordingly, 
$$
U\sim c_{f_\circ} W_{f_\circ }+ c_{f_{\partial^B}}\left( c_{f_{\partial^I},
f_{\partial^B}}W_{f_{\partial^I}}+W_{f_{\partial^B}}\right)\, ,
$$
where $C=(c_{f_{\partial^I},f_{\partial^B}})$.  
The components of the vector $\left( \begin{array}{c}
        c_{f_\circ} \\ c_{f_{\partial^B}}
        \end{array}
\right)$
are found by solving the system of second order differential equation that 
results from the weak formulation of (\ref{vp}) in this context, and which is,
of course, very closely related to (\ref{e24}). Indeed, if we carry out the 
additional splitting of the blocks of the matrices in (\ref{e24}),      
$$
\begin{array}{rcl}
I_{f_\circ,f'_{\partial}}  & = &   (I_{f_\circ,f'_{\partial^I}} \quad
I_{f_\circ,f'_{\partial^B}}) \hspace{1mm} = \hspace{1mm} \left( \begin{array}{c}
I_{f_{\partial^I},f'_\circ}  \\ I_{f_{\partial^B},f'_\circ}
\end{array}
\right)^t \hspace{2mm}  = \hspace{2mm}     
I^t_{f_{\partial},f'_\circ} \, , 
\vspace{1mm} \\
I_{f_{\partial},f'_\partial } & = & \left( \begin{array}{cc}
 I_{f_{\partial^I},f'_{\partial^I} } &
I_{f_{\partial^I},f'_{\partial^B}} \\
I_{f_{\partial^B},f'_{\partial^I}} & I_{f_{\partial^B},f'_{\partial^B}}
\end{array}
\right) = \left( \begin{array}{cc}
 I_{f_{\partial^I},f'_{\partial^I} } &
 0  \\
 0  & I_{f_{\partial^B},f'_{\partial^B}}
\end{array}
\right)
\vspace{1mm} \, , \\ 
K_{f_\circ,f'_{\partial}}  & = &   (K_{f_\circ,f'_{\partial^I}} \quad 
K_{f_\circ,f'_{\partial^B}}) \hspace{1mm} = \hspace{1mm} \left( \begin{array}{c}
K_{f_{\partial^I},f'_\circ}  \\ K_{f_{\partial^B},f'_\circ}
\end{array}
\right)^t \hspace{2mm}  = \hspace{2mm}    
K^t_{f_{\partial},f'_\circ} \, ,   
\vspace{1mm} \\
K_{f_{\partial},f'_\partial } & = & \left( \begin{array}{cc}
 K_{f_{\partial^I},f'_{\partial^I} } &  
K_{f_{\partial^I},f'_{\partial^B}} \\
K_{f_{\partial^B},f'_{\partial^I}} & K_{f_{\partial^B},f'_{\partial^B}} 
\end{array}
\right) = \left( \begin{array}{cc}
 K_{f_{\partial^I},f'_{\partial^I} } &  
 0  \\
 0  & K_{f_{\partial^B},f'_{\partial^B}} 
\end{array}
\right) 
\, ,    
\end{array}
$$
we now need to solve the Cauchy problem for the square system 
\begin{equation} \label{e29}
\rho \,  I_{\circ, \partial^B }  \left( \begin{array}{c}
\ddot{c}_{f'_\circ} \\ \ddot{c}_{f'_{\partial^B}}
\end{array}
\right)= K_{\circ, \partial^B}
\left( \begin{array}{c}
c_{f'_\circ} \\ c_{f'_{\partial^B}}
\end{array}
\right)+
\left( \begin{array}{c}
\< F_\omega, W_{f_\circ}\>  \\ \< F_\omega, c_{f_{\partial^I},f_{\partial^B}}
W_{f_{\partial^I}}+ W_{f_{\partial^B}} \>
\end{array}
\right)\, ,
\end{equation}
where
$$
\begin{array}{rcl}
I_{\circ,\partial^B} & = & \left( \! \! \begin{array}{cc}
I_{{f_\circ},{f'_\circ}} &  I_{{f_\circ},{f'_{\partial^I}}}C+
I_{{f_\circ},{f'_{\partial^B}}} \\
C^t I_{f_{\partial^I},{f'_{\circ }}}+I_{f_{\partial^B},{f'_{\circ}}} & 
C^t I_{{f_{\partial^I}},{f'_{\partial^I}}}C+I_{{f_{\partial^B}},
{f'_{\partial^B} }}
\end{array}
\! \! \right) \, , \vspace{1mm} \\
K_{\circ, \partial^B} & = & \left( \begin{array}{cc}
K_{f_\circ,f'_\circ} &  K_{f_\circ,f'_{\partial^I}}C+
K_{f_\circ,f'_{\partial^B}}  \\
C^t K_{f_{\partial^I},f'_\circ} +K_{f_{\partial^B},f'_\circ} & 
C^t K_{f_{\partial^I},f'_{\partial^I} }C+
K_{f_{\partial^B},f'_{\partial^B}} 
\end{array}
\right) \, . 
\end{array}
$$

The bottom right blocks of the matrices in (\ref{e29}) are not diagonal, as 
was so in the case for the system (\ref{e24}).
In addition to the diagonal term, these blocks contain at most three 
nonzero columns per row. This 
system is different from the (nonsymmetric) system that 
results when (\ref{e24}) is solved for 
$\left( \begin{array}{c}
        c_{f_\circ} \\ c_{f_{\partial}}
        \end{array}
\right)$
under the assumption that
$c_{f_{\partial}}=  
\left( \begin{array}{c}
        c_{f_{\partial^I}} \\ c_{f_{\partial^B}}
        \end{array}
\right)=  
\left( \begin{array}{c}
        Cc_{f_{\partial^B}} \\ c_{f_{\partial^B}}
        \end{array}
\right)$. 
Discretizing over
${\rm Div}(K)$, and then projecting the solution of the resulting 
system onto ${\rm Div}^{\rm b}(K)$,
or discretizing over ${\rm Div}^{\rm b}(K)$, and then solving the resulting
system, are not commutative operations. 

The solution of (\ref{vp}) is now approximated numerically by the vector field
\begin{equation} \label{e34}
U= c_{f_\circ} W_{f_\circ }+ c_{f_{\partial^B}}\left( c_{f_{\partial^I},
f_{\partial^B}}W_{f_{\partial^I}}+W_{f_{\partial^B}}\right)=
U= c_{f_\circ} W_{f_\circ}+ c_{f_{\partial^I}}W_{f_{\partial^I}}+
c_{f_{\partial^B}} W_{f_{\partial^B}} \, ,
\end{equation}
where the coefficients
$\left( \begin{array}{c}
        c_{f_\circ} \\ c_{f_{\partial^B}}
        \end{array}
\right)$
are given by the solution of (\ref{e29}), and we have $c_{f_{\partial^I}}=
c_{f_{\partial^I}, f_{\partial^B}} c_{f_{\partial^B}}=Cc_{f_{\partial^B}}$. 
The {\it fine} approximations to the vibration patterns of the body are  
those waves induced by the eigenvectors of the matrix  
$-(\rho I_{\circ, \partial^B})^{-1} K_{\circ, \partial^B}$, all of which
satisfy the boundary condition of (\ref{nblg}) by construction, and which
therefore correspond to waves that yield curves in ${\rm Div}^{\rm b}(K)$.
Proceeding as before, we consider a few of the negative eigenvalues of  
$-(\rho I_{\circ, \partial^B})^{-1} K_{\circ, \partial^B}$, 
and describe the resonance vibration patterns of 
(\ref{e29}) by depicting their nodal points over the portion of 
the boundary opposite to that where the external wave hits it.
We refer to these as the {\it fine resonance waves} of the body.

\section{Results} \label{s5}
In our numerical experiments, we use bodies with two types of geometries.
In each case, we describe the resonance solutions of (\ref{vp}) by the
algorithms of \S \ref{s4}, and do this for the set of frequencies 
associated to what are considered today the most important modes in tuning
violin plates.  

\begin{enumerate}
\item For the first of our experiments, $\Omega$ is a slab of  
$\text{$10$ cm}\times \text{$1$ cm} \times \text{$20$ cm}$, see Fig. \ref{f0}. 
We shall consider also thinner versions of it, of width $0.5$ cm and $0.25$ 
cm, respectively. By the flatness of the boundary, a sine-wave directed towards
the slab traveling parallel to its thin axis hits the boundary with the same
phase at all points. These simulations are used to illustrate the effect that 
this fact has on the vibration pattern, as well as the effect of rescaling the 
thin direction.

\begin{figure}[H]
\includegraphics[height=3.75in]{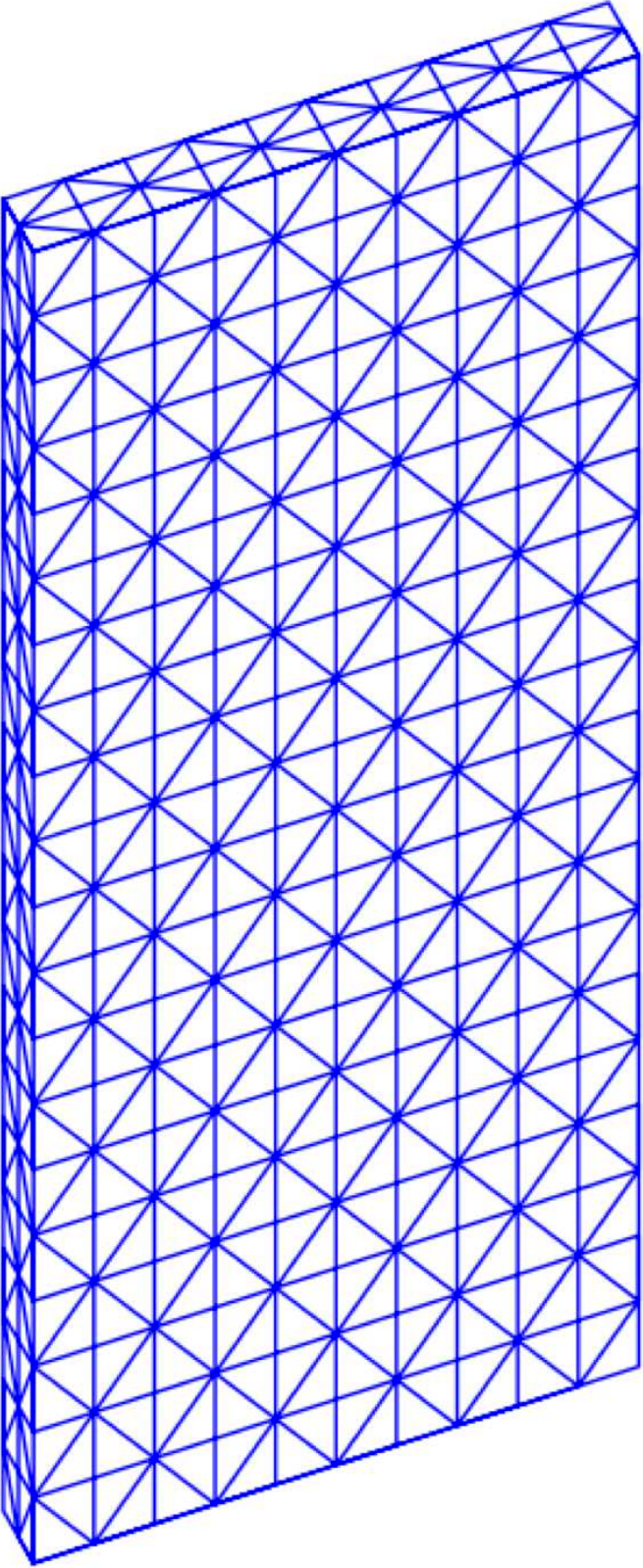} 
\vspace{1mm}
\caption{Slab with its (visible) triangulation, as used here.}\label{f0}
\end{figure}

\item For the second of our experiments, $\Omega$ has the geometry of the top 
plate of a classic violin, the Viotti, as per \cite{tstrad}. A computer 
generated view of this body from above is displayed in Fig. \ref{f1}. The
view, which does not include the $f$-holes, is 
based on measurements of the actual size images in \cite{tstrad}.   
\end{enumerate}

\begin{figure}[H]
\includegraphics[width=2.95in,angle=-90]{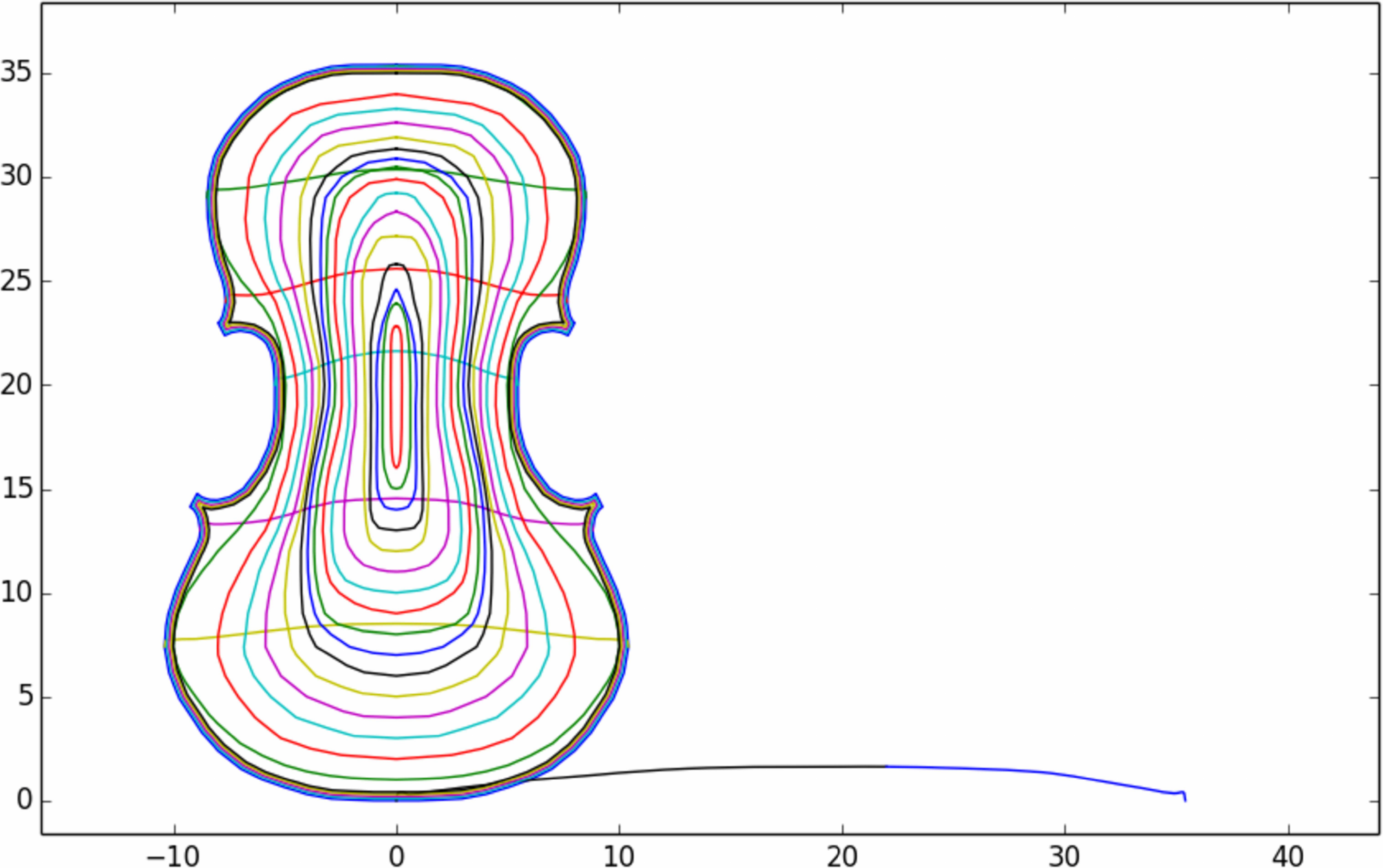}  
\caption{Viotti top plate view (without the $f$ holes) and elevation profile
curve.} \label{f1}
\end{figure}

\subsection{Computational complexity}
We subdivide the slab of 
$\text{$10$ cm}\times \text{$1$ cm} \times \text{$20$ cm}$ into $400$ regular 
blocks of size $\text{$1$ cm}\times \text{$0.5$ cm} \times \text{$1$ cm}$ each, 
Fig. \ref{f0}. The resulting triangulation $K$, where the blocks are given the
standard subdivision 
into five tetrahedrons each, contains $693$ vertices,
$3,212$ edges, $4,520$ faces and $2,000$ tetrahedrons, and
$1,040$ of the faces are boundary faces. The first barycentric subdivision
$K'$ contains $10,425$ vertices, $61,544$ edges, $99,120$ faces (of which
$6,240$ are boundary faces), and $48,000$ tetrahedrons.
The two other thinner versions of the slab have triangulations with the same 
number of elements, merely scaling the depth accordingly. 
The aspect ratio of their tetrahedrons are half and a quarter of the aspect
ratio of those in the first of the triangulations, respectively. 

When $\Omega$ is the top plate of the Viotti violin, the triangulation that
our work requires is significantly more complex. The top view in Fig. \ref{f1}, 
which depicts the outline of this plate, 
includes the purfling curves, the transversals A, B, C, D, and E, 
sixteen level sets, and the elevation curve that combines the F and G curves, 
as per  measurements of the images in \cite{tstrad}.   
The elevation curve serves as a reference for the height of the 16 
level sets in the view, with the level of the $j$th set being the elevation 
at a distance of $j$ cm along it, starting at $0$.   
This $\Omega$ can be inscribed into a rectangular box of
$35.4\,  {\rm cm} \times 20.8 \,  {\rm cm}$. This body is curved, as
the 5 transversals to the level sets in the Fig. suggest. The thickness of the
plate varies nonuniformly, ranging from a lowest of $0.21$ cm to a largest of 
$0.36$ cm.    

In order to get a reasonable resolution for the type of thickness and
curvature that this body has, away from the edge, we subdivide 
$\Omega$ into blocks of
size $\text{$0.5$ cm}\times \text{$0.5d$ cm} \times \text{$0.5$ cm}$,
where $d$ is an average thickness of the plate at points where the block is
located. We proceed similarly at the edge, but the first dimension of the 
blocks we consider there is taken to be nonuniform, out of necessity.
 All in all, it takes 4,608 of these blocks to cover $\Omega$, and
the triangulation $K$ that we then derive contains 7,170 vertexes, 
35,157 edges, 51,028 faces (of which 10,232 are boundary faces), and 
23,040 tetrahedrons. The first barycentric subdivision of this 
triangulation has 116,395 vertexes, 699,042 edges, 1,135,608 faces, and 
552,960 tetrahedrons, with 61,392 of the faces located on the boundary.

In either case, as we traverse the body across the two blocks separating 
the bounding surface in the thin direction, the barycentric triangulations
used contain 444 faces, with 35 intermediate faces separating pair of 
boundary faces on opposite ends, which when included yield a total 37 faces 
altogether in going from one side to the other. This provides an adequate 
resolution for the numerical approximation to accurately capture the true 
nature of the vibrating wave in these regions, and to propagate
in all transversal directions it goes. Refinements of the triangulations 
that we use would improve the accuracy of the numerical solution, but such 
would lead to a complexity that is out of the scope of our current 
computational resources.

As for the boundary conditions that are satisfied by the fine 
vibration waves approximations, of all the $2 \times 6$ subsystems of
(\ref{e27}) that the faces in ${\mc F}_{\partial }$ generate, 
of the 1,040 of them for the slab, exactly eleven 
have rank one, and of the 10,232 of them for the Viotti plate, exactly
one have this property. The differences between these numbers reflect the
$\mb{Z}/2 \times \mb{Z}/2 \times \mb{Z}/2$ symmetry of the triangulation
of the slab, as opposed to the nonsymmetric nature of the triangulation of
the violin plate, whose thickness varies in a nonuniform manner throughout
the body. The astute reader should have noticed that the exterior normal to a
boundary face, and its normal as an oriented two simplex, are not 
necessarily the same.

We find the matrices $I$ and $K$ of the systems (\ref{e24}), (\ref{e29}),  
by executing python code written for the purpose. The code is structured into
four major modules, where just the last one is body specific, an
object-oriented conception aimed at making it easy to expand our analysis 
to bodies with geometries other than those considered here. 
The processing of the results, and their graphical display, is carried out
by some additional python code written for the purpose, the graphical 
component of it built on top of the PyLab standard library module.

The splitting into interior and boundary faces leads naturally to generating 
separately the corresponding blocks of these matrices.
With a 2.4GHz Intel Core i7 processor, and 8GB 1600 MHz of memory, generating
the upper left blocks of the matrices for the slabs takes close to 7h of
CPU time each, about 36m for the upper right, and close to 2s for the
bottom right. By contrast, for the violin plate, in 10h43m of CPU time
we generate 0.59\% of the upper left blocks, while in 20h30m we generate
29.35\% of the upper right block. The bottom right blocks of the
matrices take a time comparable to their alteregos for the slabs. These times
are not optimal; they can be improved using the ideas behind the 
parallelization of the calculations on the basis of the linear ordering of 
the faces. This, in fact, was the way we managed to complete the generation
of the said matrices in the said hardware.

\subsection{Elastic constants}
We assume that the material that makes our bodies is orthotropic with the 
elastic constants of spruce. Notice that the arcs of the wedges of 
either maple or spruce used for a two-piece violin plate subtend a very 
small angle, and so the transformation from cylindrical to Cartesian 
coordinates on these wedges is given by a matrix that is very close to the 
identity, see Remark \ref{r6} here. Consequently, we take the components 
of tensors expressed in these two coordinate systems to be the same.   

The elastic constants values we use are those  in tables 1 \& 2 for the 
Engelmann Spruce. Since these values reflect a failing condition 
(\ref{rel}), we 
take the average of the computed values of
$\mu_{ij}/e_i$ and $\mu_{ji}/e_j$ as the value of either one of 
these quantities in our calculations. Thus, the diagonal blocks
of the tensor $W$ of elastic constants (\ref{ten}) in our simulations are 
$$
W_{3\times 3}=10^7 \left( \begin{array}{ccc}
157.198269069862 & 44.1920517114940 &  116.065341927474  \\ 
44.1920517114940 & 72.0200103705017 & 75.6887031695923  \\ 
116.065341927474 & 75.6887031695923 & 1095.80735919001  
\end{array}\right) \, ,
$$ 
and
$$
D_{3\times 3} =10^7 \left( \begin{array}{ccc}
117.480 & 0 & 0 \\ 
0 & 121.396 & 0 \\ 
0 & 0 & 9.790 
\end{array}\right) \, ,
$$
respectively, the units of the components in Pa.  
Since we are using  
an orthonormal frame, we can raise or lower indices in tensors with abandon.
We use $\rho = 360 \; {\rm kg}/{\rm m}^3$ for the density parameter. 
Notice that the eigenvalues of $W_{3\times 3}$ are 
$$
10^7\{ 156.292790395160, 1116.15681336097, 52.5760348742398 \} \, , 
$$
so the stored energy function of each of our bodies is coercive.

\subsection{Simulations} 
For each of the bodies with their triangulations as above, we analyze the 
divergence of the eigenvector solutions, and resonance waves, of the 
systems (\ref{e24}) and (\ref{e29}), respectively.

\subsubsection{The divergence of the coarse and fine normalized 
eigenvector solutions} 
The eigenvalues and eigenvectors of the homogeneous systems
associated to (\ref{e24}) and (\ref{e29}) are generated using the 
ARPACK routine \verb+eigsh+ in shift-invert mode, with their corresponding 
matrix parameters $\rho I$, and $K$, respectively.  This computes the 
solutions $(\lambda, c)$ of the system 
$$
\rho I c = \lambda K c \, .
$$
With \verb+sigma=-1/(2 pi f)^2+, and \verb+which='LM'+ passed onto 
\verb+eigsh+, we execute the routine for a frequency 
\verb+f+ any of 
$80$, $147$, $222$, $304$, and $349$ Hz, respectively. In each case, we 
produce pairs \verb+-((2 \pi f_r)^2, c^{f_r})+ of eigenvalue and 
eigenvector for the matrix $(\rho I)^{-1}K$, where \verb+f_r+ is the eigenvalue
of the matrix that is closest to the inputted \verb+f+, in magnitude, and let 
$$
U_{f_r}=e^{2\pi f_r t i}\sum_{W_f} c_f W_f=e^{\omega_{f_r}t i}\sum_{W_f} c_f W_f
$$ 
be the corresponding eigenvector wave solution. 
In the case of the system (\ref{e24}), we let 
$$
U_{f^f_{coarse}}=e^{2\pi f_{coarse} t i}\sum_{W_f\in \rm{Div}(K)}
c^{f_{coarse}}_f W_f
$$
be the {\it normalized} coarse eigenvector solution, while in the case of the
system (\ref{e24}), we let 
$$
U_{f^f_{fine}}=e^{2\pi f_{fine} t i}\sum_{W_f\in {\rm Div}^b(K)} 
c^{f_{fine}}_f W_f
$$
be the {\it normalized} fine eigenvector solution that 
the said pair produces.

We study the divergence of any normalized eigenvector solution 
$U=e^{2\pi f_r t i}\sum_{W_f} c_f W_f$ by computing the flux 
$$
\int_{\partial \Omega} U \cdot n \, d\sigma 
$$
through the boundary of the body at time $t=0$. For the 
coarse and fine eigenvector solutions, at the frequencies above and for
each of the bodies under consideration, the results are as follows:
\medskip

\begin{center}
\begin{tabular}{|c|r|r|r|r|r|}
\hline \hline
Body & $f$\phantom{3} & $f_{coarse}$\phantom{333} & flux $U_{f^f_{coarse}}$ & $f_{fine}$\phantom{333} & flux $U_{f^f_{fine}}$ \\ \hline \hline
\multirow{5}{*}{Slab 1.0} & 80 & 79.89682695 & 0.0066641031 & 79.89486045 & 0.0012710885 \\
& 147 & 146.81041954 & -0.0066636678 & 146.81041954 &-0.0012710492 \\
& 222 & 221.71369483 &  0.0066635620 & 221.71369483 &-0.0012710349 \\
& 304 & 303.60794252 &  0.0066635192 & 303.60794252 &-0.0012710230 \\
& 349 & 348.54990773 &  0.0066635053 & 348.54990774 &-0.0012710162 \\ \hline
\multirow{5}{*}{Slab 0.5} & 80 & 79.89682695 & 0.0062420055 & 79.89682695 & 0.0032331563  \\
 & 147 & 146.81041954 &-0.0062419805 & 146.80677110 & 0.0032331603 \\ 
 & 222 & 221.71369483 & 0.0062419720 & 221.71369483 &-0.0032331560 \\ 
 & 304 & 303.60794251 & 0.0062419656 & 303.60794251 &-0.0032331470 \\ 
 & 349 & 348.54990771 & 0.0062419620 & 348.54990772 &-0.0032331407 \\ \hline  
\multirow{5}{*}{Slab 0.25} & 80 & 79.89682695 & 0.0047848418 & 79.89682695 & -0.0044987163 \\
 & 147 & 146.80681943 & 0.0047851684 & 146.81041954 & -0.0044986653 \\
 & 222 & 221.71369483 & 0.0047852434 & 221.71369483 & -0.0044986514 \\
 & 304 & 303.60794250 & 0.0047852684 & 303.60794251 & -0.0044986436 \\
 & 349 & 348.54179550 & 0.0047852740 & 348.54990771 & -0.0044986401 \\ \hline
\multirow{5}{*}{Viotti plate} & 80 & 79.89682695 & -0.0026727377 & 79.93578386 & 0.0350092360 \\
 & 147 & 146.81041953 & 0.0026739488 & 146.99894842 &  0.0350815837  \\
 & 222 & 221.71369480 & 0.0026742324 & 221.81928975 & -0.0033142755 \\
 & 304 & 303.60794243 & 0.0026743267 & 304.01369231 &  0.0176140244 \\
 & 349 & 348.54990760 & 0.0026743492 & 348.54990764 & -0.0064094032 \\
\hline \hline
\end{tabular}
\smallskip

\centerline{Table 3. Initial fluxes for the coarse and fine eigenvector 
solutions.}
\end{center}

\subsubsection{Coarse and fine resonance waves} 
We subject the body to an external sinusoidal pressure wave of the form 
$$
\vec{F}=\vec{F}_0 \sin({\bf k}\cdot {\bf x} \mp \omega t) 
=\vec{F}_0 \sin({\bf k}\cdot {\bf x})\cos({\omega t})\mp
\vec{F}_0 \cos({\bf k}\cdot {\bf x})\sin({\omega t}) 
$$
that travels in the appropriate direction for it to hit the bottom 
of the slab, or belly of the Viotti plate, first. If $k=|{\bf k}|$ is the
magnitude of the wave vector, we use $v=\omega/k=343$ m/sec,  
the speed of sound in dry air at $20\mbox{}^\circ$C. The source of the
wave is placed at a distance of $62$cm from the body, along a line 
that passes through the height-width plane of the body perpendicularly at 
the half way point of both, its height and width. This external force 
induces a force $\vec{F}_{\omega}$ on the body.  

The coarse (\ref{e24}) and fine (\ref{e29}) systems are considered with the  
nonhomogeneous force term arising from the $\vec{F}_{\omega}$ on the body
induced by the external wave $\vec{F}$ above, and with trivial initial data. 
(These are the coarse and fine discrete versions of the initial value problem 
for (\ref{vp}).)  The nonhomogeneous terms in these systems are vectors
of the form
$$
\vec{C}_1 \cos{(\omega t)}\mp \vec{C}_2 \sin{(\omega t)}\, ,
$$
where $\vec{C}_1$ and $\vec{C}_2$ are time independent vector fields on the
body. 
If $\omega =2\pi f$ for $f$ any of the frequency values 
$80$, $147$, $222$, $304$, and $349$ Hz, respectively, 
we let $\omega_{f_r} =
2\pi f_r$ be the closest eigenvalue to $\omega$ that the matrix 
$(\rho I)^{-1}K$  
has, $\rho I$ and $K$ the matrices of the system in consideration.
The resonance wave that this vibration mode produces is given by  
\begin{equation} \label{resw}
W_{f_r}=\frac{1}{\omega^2_{f_r}-\omega^2}\sum_{W_f}
\left(c^{W_f}_1\left(-\cos{(\omega_{f_r}t)}+
\cos{(\omega t)}\right)\pm c^{W_f}_2 \left( \frac{\omega}{\omega_{f_r}}
\sin{(\omega_{f_r}t)}-\sin{(\omega t)}\right)\right)W_f \, ,
\end{equation}
where, for each $j=1,2$, the vector $\vec{c}_j =(c_j^{W_f})$ is the solution 
to the linear system of equations  
$$
-((2\pi f)^2 \rho I +K)\vec{c}_j = \vec{C}_j \, . 
$$
The summation is over the basis elements $W_f$ of ${\rm Div}(K)$ and
${\rm Div}^b(K)$ for the coarse and fine systems, respectively.
We generate the vector $\vec{C}_j$ as a function of $\vec{F}$, and the geometry
of the body, and then solve the system of equations above for $\vec{c}_j$ 
using the \verb+scipy.sparse.linalg+ routine \verb+spsolve+, with the 
appropriate parameters.

We compute the values of the resonance wave at the barycenter, and vertices, 
of the boundary faces on the boundary side of the body opposite to the 
incoming external wave.
There are $3,661$ such points for the slabs, and 
$42,069$ for the Viotti plate. We do these calculations at the equally spaced 
times $t_j=j\frac{2\pi}{10\omega}$, $j=1 , \ldots, 10$, corresponding to a 
full cycle of the external wave. For each $t_j$, we find the maximum
${\rm max}_{t_j}=\max{\{\| W_{f_r}(t_j)\|}\}$ and minimum 
${\rm min}_{t_j}=\min{\{\| W_{f_r}(t_j)}\|\}$ of the set of norms of the 
resonance wave at the indicated points, and with 
$\delta_{t_j}=\frac{1}{10}({\rm max}_{t_j}-{\rm min}_{t_j})$,  
any of the said points is considered to be nodal at time $t_j$ if the norm 
of the resonance wave solution at the point is no larger than 
${\rm min}_{t_j}+ c_\Omega \delta_{t_j}$, with $c_\Omega=0.8$ for the slabs,
and $c_\Omega=0.04$ for the Viotti plate, respectively. A point is 
defined to be nodal if it is nodal at all the $t_j$s. 
All the results for the coarse and fine resonance waves are depicted
in Figs. 3-6 below. In each case, we indicate
the value of $c_\Omega$ that is being used to define a point as nodal.    

We study the change in the resulting resonance pattern produced by taking
into consideration the six modes with 
eigenvalues $\omega_{f_{r_j}}=2\pi f_{r_j}$, $j=1,\ldots, 6$, closest to 
$\omega=2\pi f$, as opposed to the single closest one, as above. The 
resonance wave is then a sum of six terms as 
in (\ref{resw}), one per $f_{r_j}$s. We exhibit the results for 
$f=147$Hz in Fig. \ref{f7}. The five cases depicted correspond to nodal
points defined as above for values of 
$c_\Omega=0.04, 0.02, 0.01, 0.005$, $0.0025$,
respectively. By comparing the corresponding cases in Figs. 6 and 7, we 
observe now a better definition of the details of the coarse and fine 
resonance patterns, though the patterns themselves are not changed 
significantly (the actual number of nodal points in each of these cases came 
out to be exactly the same, surely, a coincidence). On the other hand,
as the notion of a nodal point becomes stricter by decreasing the value of
$c_\Omega$, clear details in the resulting patterns emerge, and appear  
to point quite closely towards the image by holographic interferometry of
the mode 2 of a top violin plate in \cite[Fig. on p. 177]{hu}. 
Since of the frequencies we consider, the one at which our model for the 
Viotti plate vibrates the poorest is 147Hz, in spite of the differing 
conditions between our simulations and the experiments in \cite{hu}, we 
take the favorable comparison just made as a validation of our results. 
The comparison gets better for any of the other values of $f$.

In all of our simulations, for a given geometry and for a given frequency,
the fine resonance waves have fewer nodal points than the coarse. For the 
slabs, the number of nodal points of the resonance waves decreases as the 
frequency increases, with their values and changes from one to the
next smaller and smaller as the thickness decreases.  By contrast, 
these last results are markedly different for the Viotti plate, for which 
the number of nodal points of the coarse resonance waves decreases as the 
frequency increases, but the changes are not monotonically decreasing, 
and the largest of them occurs in going from 147Hz to 222Hz. The
pattern is broken altogether for the fine resonance waves, for which, by far, 
the largest number of nodal points occurs at 147Hz; this number decreases then 
monotonically for 222Hz and 304Hz to values less than the value it has at 80Hz, 
and then goes up slightly for the wave associated to 349Hz. The phenomena
appears to be more than merely an issue of curvature, perhaps related to 
wavelength and thickness also, or the combination of (at least) these 
three elements.  
  
For our simulations for the slabs indicate a limitation to the effect 
of selective thinning of the plate, a practice among luthiers that is 
commonly believed to lower 
or increase frequency modes if carried out in regions of high or low 
curvature, respectively. The values of the ten quotients 
${\rm max}_{t_j}/{\rm min}_{t_j}$, 
ordered from smallest to largest, align with thicknesses of 1cm, 0.25cm, 
and 0.5cm, respectively. At some point, the effect of additional thinning seems 
to be reversed, and though it would be hard to test this reversal on the 
already quite thin violin plates, if at all possible, the fact that the 
Viotti plate has large portions of it that are thinner than the thinnest of the 
slabs, and the nonmonotonic property of the number of nodal points of its 
resonance waves, is consistent with this result for the slabs.
For the slabs, the largest of the alluded quotients, at a given one of those 
frequencies used here, is in the order of $10^3$; for the Viotti plate, it 
is in the order of $10^{11}$. 

For lack of data on the elastic constants of actual wood used in the making
of violin plates, we have not studied possible changes in the results due to 
the relative humidity of the material. If compared to the original, our model 
for the Viotti plate should be rather dull, and vibrate poorly, as the data
we use is from wood with 
a high 12\% moisture content. 
(A renowned luthier has told us that the wood he uses to make the 
plates has moisture content in the range of 2\%-5\% \cite{ca}. This content 
is likely lowered when the plate is coated with varnish, which adds mass and 
absorbs some moisture as it dries. Note that the addition of the varnish 
changes also the flexibility by stiffening the plate across the grain.) 
We could analyse also temperature increases in the body due to 
dissipation, or friction, 
and how they affect the vibration patterns of the plate, but the results of 
such analyses are out of the scope of the article, as is the study of the
resonance pattern under the influence of large exterior pressure forces.

We denote by $W_{f_{coarse}^f}$ and $W_{f_{fine}^f}$ the normalized coarse
and fine resonance wave solutions associated to the pair $(f,f_r)$ for
(\ref{e24}), and (\ref{e29}), respectively. As the waves
start with trivial initial condition, we quantify the extent to
which our algorithms maintain the divergence free condition throughout time
by  evaluating their fluxes over the boundary at the time $t_j$ where
the quotient ${\rm max}_{t_j}/{\rm min}_{t_j}$ is the largest, and which
happens to be $t_7$ in all of our simulations. The
normalization performed on the resonance waves
make these results independent of the magnitude of $\vec{F}_0$
in the external wave $\vec{F}$ that induces the resonance. They are
listed in Table 4.
\medskip

\begin{center}
\begin{tabular}{|c|r|r|r|r|r|}
\hline \hline
Body & $f$\phantom{3} & $f_{coarse}$\phantom{333} & flux $W_{f^f_{coarse}}$ & $f_{fine}$\phantom{333} & flux $W_{f^f_{fine}}$ \\ \hline \hline
\multirow{5}{*}{Slab 1.0} & 80 & 79.89682695 & 0.4115594102 & 79.89486045 & 0.2154917860 \\
& 147 & 146.81041954 & 0.4077021407 & 146.81041954 & 0.2159199992 \\
& 222 & 221.71369483 & 0.4009923846 & 221.71369483 & 0.2144771805 \\
& 304 & 303.60794252 & 0.3930235693 & 303.60794252 & 0.2119845248 \\
& 349 & 348.54990773 & 0.3887672035 & 348.54990774 & 0.2104879098 \\ \hline
\multirow{5}{*}{Slab 0.5} & 80 & 79.89682695 & 0.4264112637 & 79.89682695 & 0.2340175942 \\
 & 147 & 146.81041954 & 0.4225227485 & 146.80677110 & 0.2341400584 \\ 
 & 222 & 221.71369483 & 0.4156736460 & 221.71369483 & 0.2322834707 \\ 
 & 304 & 303.60794251 & 0.4075092711 & 303.60794251 & 0.2293410711 \\ 
 & 349 & 348.54990771 & 0.4031419661 & 348.54990772 & 0.2276115430 \\ \hline
\multirow{5}{*}{Slab 0.25} & 80 & 79.89682695 & 0.4312931255 & 79.89682695 &
0.2097749790 \\
 & 147 & 146.80681943 & 0.4273781672 & 146.81041954 & 0.2096177445 \\
 & 222 & 221.71369483 & 0.4204581119 & 221.71369483 & 0.2077423096 \\
 & 304 & 303.60794250 & 0.4122164484 & 303.60794251 & 0.2049398128 \\
 & 349 & 348.54179550 & 0.4078257909 & 348.54990771 & 0.2033199948 \\ \hline
\multirow{5}{*}{Viotti plate} & 80 & 79.89682695 & -0.1850124890 & 79.93578386 & -0.0925756811 \\
 & 147 & 146.81041953 & -0.1091897994 & 146.99894842 &  -0.0189475658 \\
 & 222 & 221.71369480 & -0.0928110618 & 221.81928975 &  -0.0426250223 \\ 
& 304 & 303.60794243 &  -0.1146617479  & 304.01369231 & -0.0544211220 \\
& 349 & 348.54990760 &  -0.1354831097 & 348.54990764 &  -0.0709637621 \\ 
\hline \hline
\end{tabular}
\smallskip

\centerline{Table 4. Fluxes of the normalized coarse and fine resonance 
waves at time $t_7$.}
\end{center}

\section{Concluding remarks}\label{s6}
If we put the computational complexities aside, now we are able to find 
numerical solutions to the Cauchy problem for the equation of motions of 
incompressible elastodynamic bodies (\ref{eq1}), (\ref{eq2}), (\ref{eq3})
themselves. We follow the proof of Theorem \ref{th2} of \cite{ebsi2}, and
to ignore (\ref{eq1}),  
choose a sufficiently large constant $\lambda$ to modify (\ref{eq2}) to
$$
\rho \ddot{\eta }(t)(x)={\rm Div}\, W^{'}(D\eta (t)(x))+
J(\eta(t))\nabla q +\lambda J(\eta(t)) \nabla_\eta J(\eta(t) ) \, ,
$$
and modify (\ref{eq3}) by adding the condition $J(\eta )=1$ on 
$\partial \Omega$. Let $(\BOne, w)$, ${\rm div}\, w =0$ be the 
Cauchy data, and view the Cauchy problem for the equation above as a 
nonlinear evolution equation
$$
\frac{d}{dt}\left( \begin{array}{c}
                   \eta \\
                   \dot{\eta}
\end{array} \right) 
= G_\lambda(\eta,\dot{\eta}) 
$$
with that initial data. If we linearize this equation at a curve
$(\eta(t),\dot{\eta}(t))$ satisfying the initial conditions, we 
obtain a quasilinear system for an unknown $u$ satisfying Cauchy 
data compatible with $(\eta(0),\dot{\eta}(0))=(\BOne, w)$. We 
may apply  the algorithm of \S \ref{s4}, conveniently changed to account 
for the modifications in the equations, and obtain a numerical solution 
$U=U(t)$ to this system. If for sufficiently large $\gamma$ we then 
consider the equation
$$
\gamma \left( \begin{array}{c}
                   \zeta \\
                   \dot{\zeta}
\end{array} \right) 
- G_\lambda(\zeta,\dot{\zeta}) = 
-U(t)+\gamma\left( \left( \begin{array}{c}
                   \zeta(0) \\
                   \dot{\zeta}(0)
\end{array} \right)
+\int_0^t U(s)ds\right)  \, , 
$$
and solve it for $(\zeta,\dot{\zeta})$ for each fixed $t$, the resulting 
mapping  $(\eta(t), \dot{\eta}(t)) \rightarrow (\zeta(t),\dot{\zeta}(t))$ 
completes the first step of the Newton scheme argument that yields the 
solution of the Cauchy problem for (\ref{eq1}), (\ref{eq2}), (\ref{eq3}). 
Iterations of the scheme would  produce a sequence $\eta_n(t)$ that 
converges to a solution $\eta$ of the problem on some time interval $[0,T]$.
Numerically, we just implemented the first step of this scheme. Should we
be able to iterate it at least one more time, we would get 
a fairly good approximate solution to the free boundary value problem under
consideration, with the diffeomorphism numerical solution being reasonably 
closed to one that preserves volume everywhere.

Any solution $u$ to the linearized equations of motion for 
(\ref{eq1}), (\ref{eq2}), (\ref{eq3}) about a curve $\eta$ is such that
${\rm div}_\eta u=0$, and so if continuous in time, it yields a
curve in the Abelian group $H^2(M;\mb{R})$. Our success in overcoming 
the computational
complexities of the problems treated here is in great part the consequence of 
the algebraic topology encoded into the Whitney forms of the triangulation.
It is quite hard to maintain a closed condition as you solve numerically
any equation. But if you discretize $u$ in spaces that are natural relative to
the cohomology class it represents, the well-posedness of the equation
will maintain your numerical solution within reasonable limits of that class 
throughout time if that condition is made to hold at the start. That 
remains so when extending the analysis to bodies that have edge and
corner singularities. The results of our simulations validate that assertion.   

The use of Whitney forms associated to edges
is well established \cite{rao, asv,boss}, though often they are employed 
to discretize physical quantities that truly represent cohomology classes in 
degree two, rather than one, and for which the use of the Withney forms 
associated to faces would be more natural instead.
Historically, functions have been approximated by computing their values 
at sufficiently many points, rather than expanding them as linear combinations
of the elements in the partition of unity given by the Whitney forms  
associated to vertices. Regardless, the use of all of these forms 
 (advocated recently by some \cite{boss2}, though Whitney
himself had used them already for several ``computational''  purposes) is    
quite natural. If their properties are exploited well, 
we may capture the essential algebraic structure underlying the problems under 
consideration, and, with a minimal number of computational elements, derive
accurate results for problems with large intrinsic complexities. 
The geometric content of a triangulation is a powerful tool to use 
to compute polytope quantities of physical significance \cite{gsi}; the power 
of this tool is several times fold larger if we include in the considerations 
its algebraic content as well.

As indicated earlier, our innovative approach can be extended to the 
study of vibrational patterns of elastic smooth bodies with corners, under 
mild assumptions on the stored energy function, 
generalized Hookian bodies for instance. This is the case of the orthotropic 
slab that we considered. If the variational principle used to derive the 
equations of motion is extended to treat cases where part of the boundary is 
fixed, then we could treat various plates vibrational problems effectively, a 
cantilever, or several others \cite{mahi, benn, bous, boud, elha}, with the 
same degree 
of generality. And since we can obtain the damping vibration modes 
as easily as the free ones while the motion remains elastic, we 
could treat the changing temperatures 
brought about by a dissipation of energy within 
the body, which in turn brings about changes in the elastic constants, 
and could eventually lead to breakages if the body is pushed all the way to
the nonlinear plastic regime. Our method allows for the use of fundamental 
principles of physics when solving these problems.

\begin{figure}[H]
\includegraphics[height=5.2in,angle=-90]{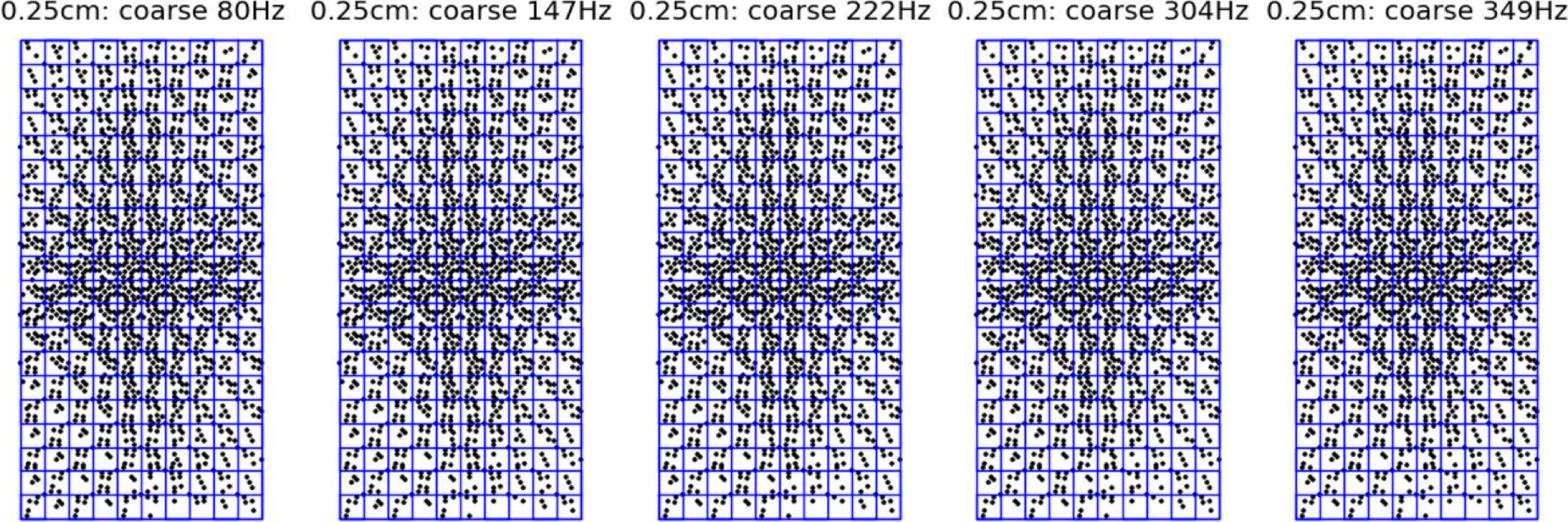}

\includegraphics[height=5.2in,angle=-90]{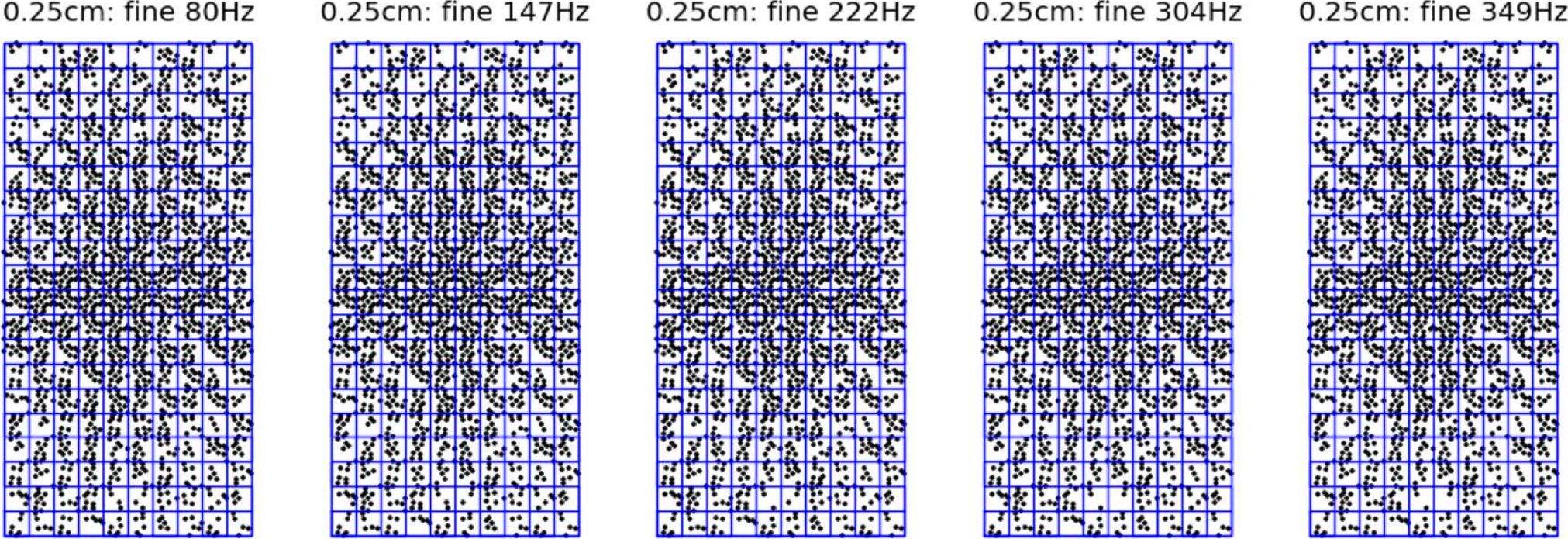}
\smallskip

\caption{Slab $10\times 0.25 \times 20$: Nodal points of the resonance waves 
at $f$Hz arising from the slab mode of vibration of frequency $f_r$ closest 
to $f$. The notion of nodal point is defined by choosing $c_\Omega=0.8$.}
\label{f3}
\end{figure}

\begin{figure}[H]
\includegraphics[height=5.2in,angle=-90]{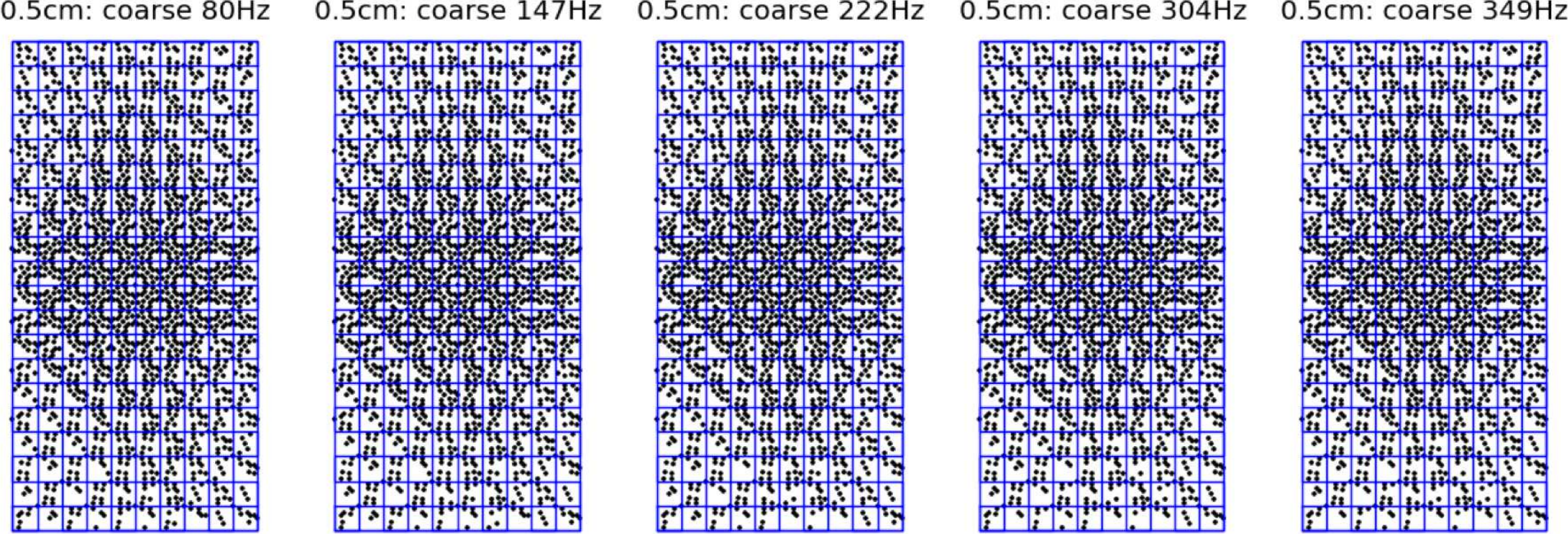}

\includegraphics[height=5.2in,angle=-90]{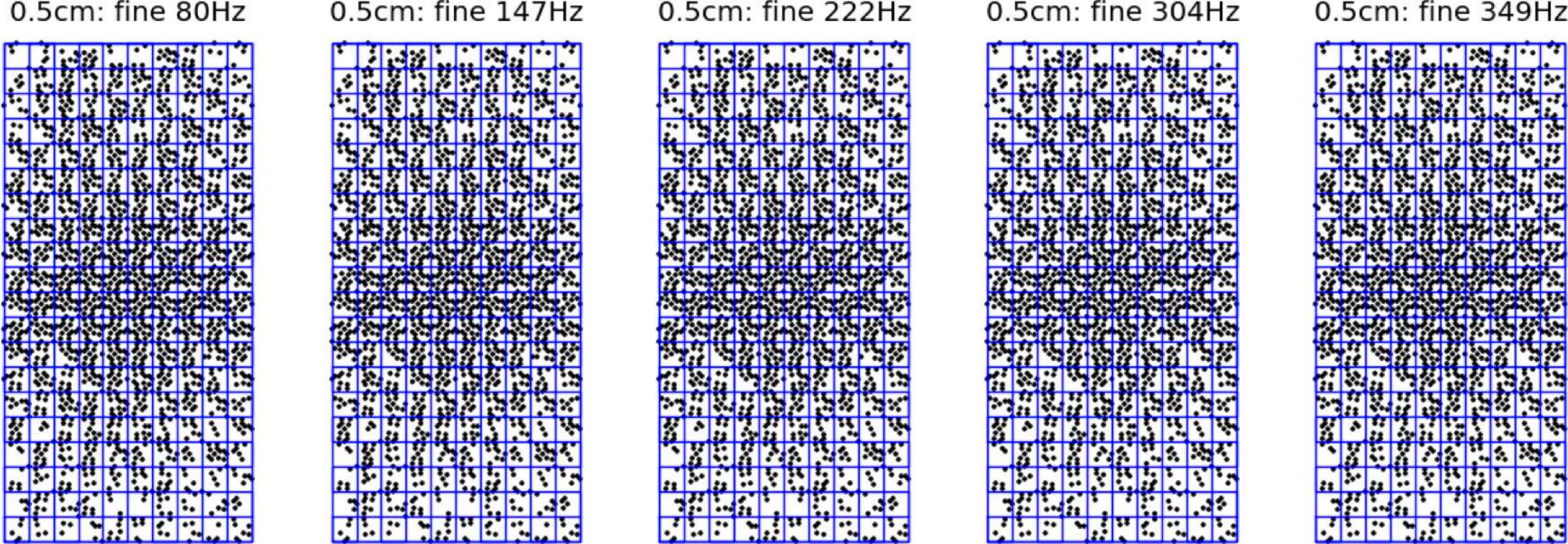}
\smallskip

\caption{Slab $10\times 0.5 \times 20$: Nodal points of the resonance waves 
at $f$Hz arising
from the slab mode of vibration of frequency $f_r$ closest to $f$. The notion of
nodal point is defined by choosing $c_\Omega=0.8$.}
\label{f4}
\end{figure}

\begin{figure}[H]
\includegraphics[height=5.2in,angle=-90]{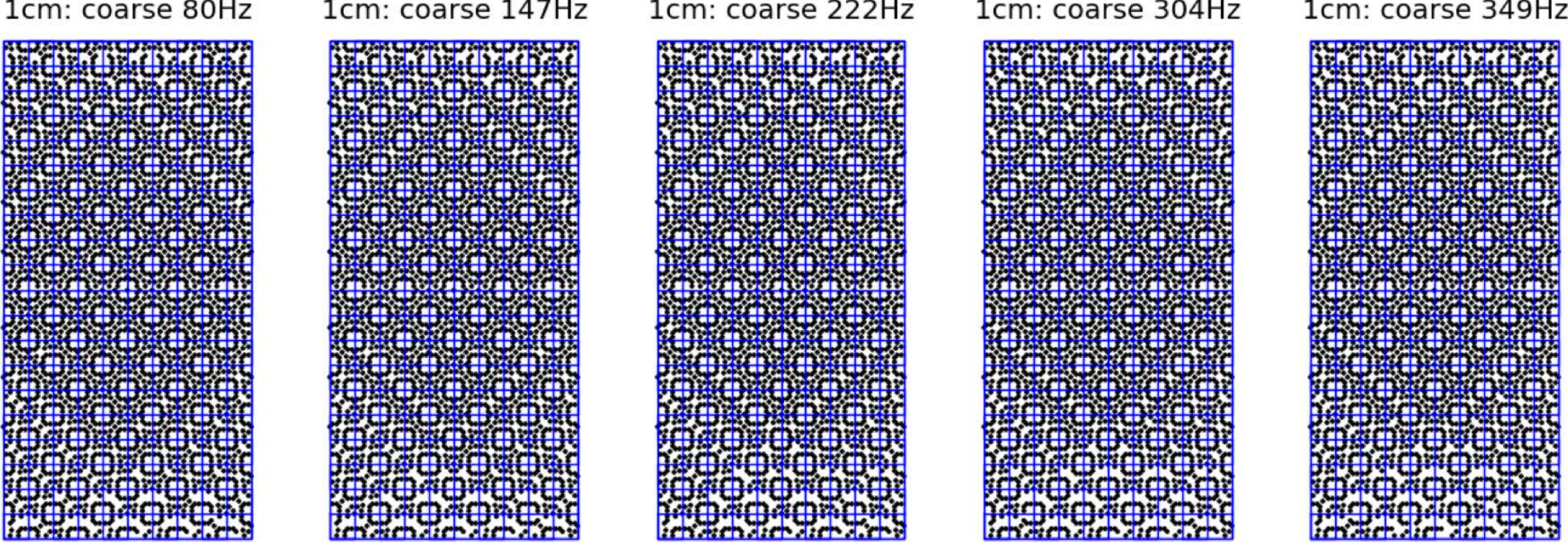}

\includegraphics[height=5.2in,angle=-90]{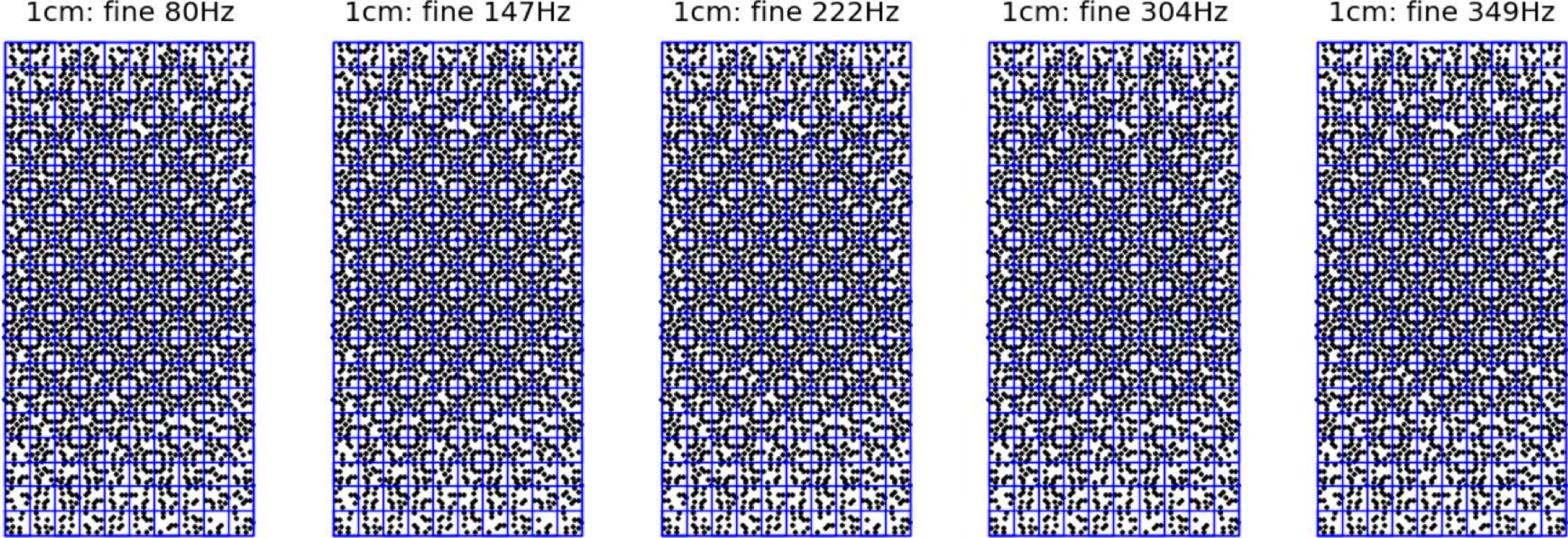}
\smallskip

\caption{Slab $10\times 1.0 \times 20$: Nodal points of the resonance waves 
at $f$Hz arising
from the slab mode of vibration of frequency $f_r$ closest to $f$. The notion of
nodal point is defined by choosing $c_\Omega=0.8$.}
\label{f5}
\end{figure}

\begin{figure}[H]
\includegraphics[height=5.2in,angle=-90]{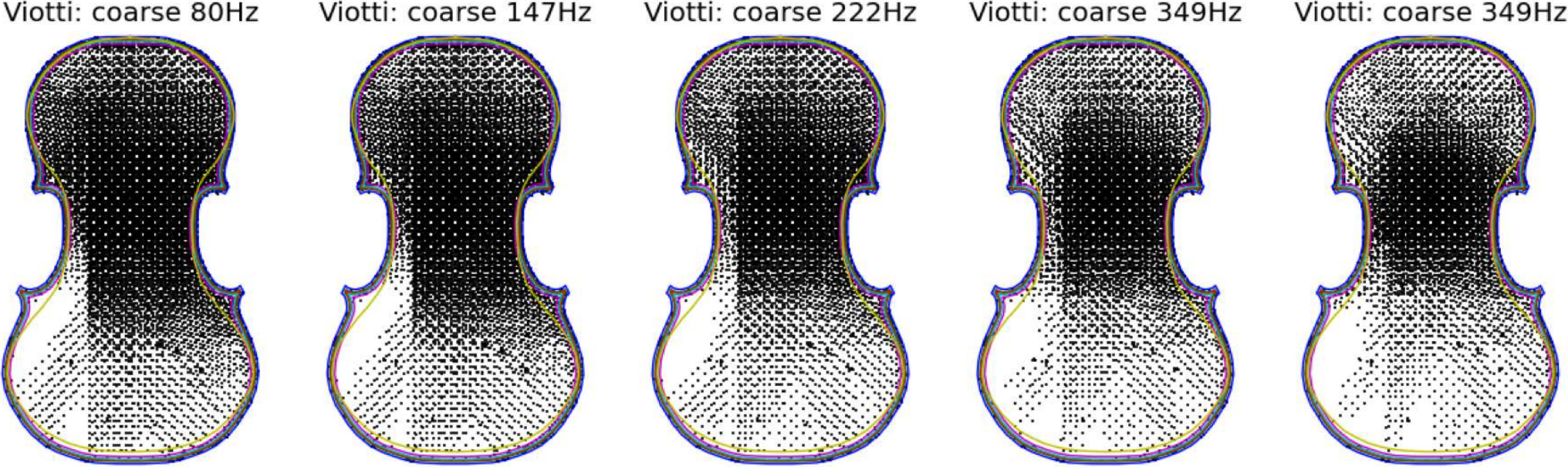}
\smallskip

\includegraphics[height=5.2in,angle=-90]{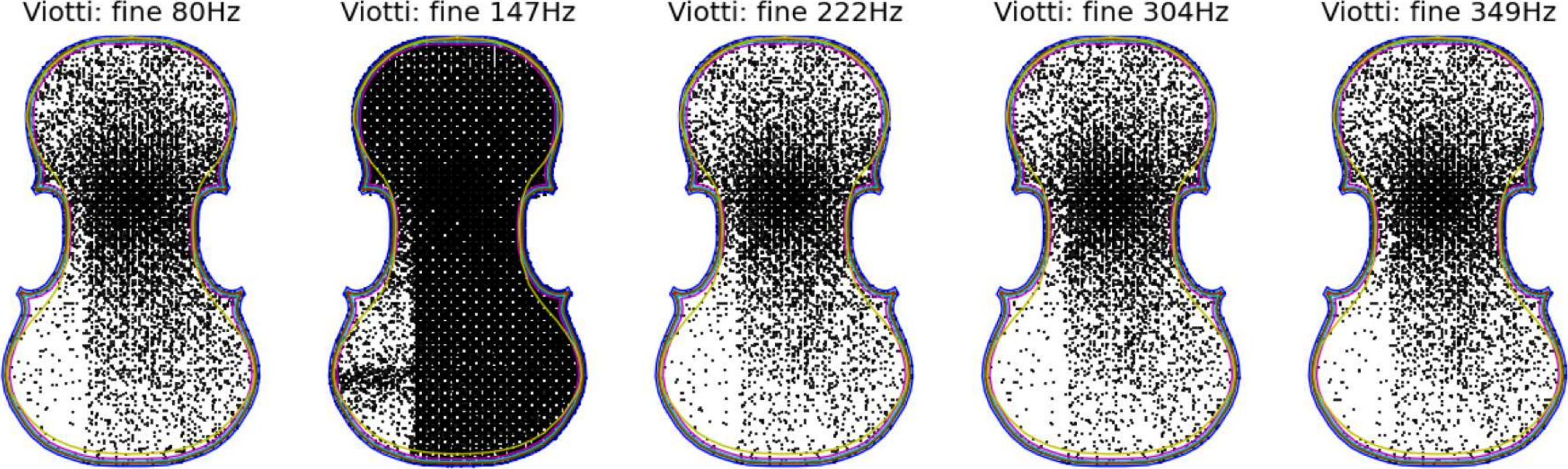}
\smallskip

\caption{Viotti plate: Nodal points of the resonance waves at $f$Hz arising 
from the mode of vibration of the plate of frequency $f_r$ closest to $f$. The 
notion of
nodal point is defined by choosing $c_\Omega=0.04$.}
\label{f6}
\end{figure}

\begin{figure}[H]
\includegraphics[height=5.2in,angle=-90]{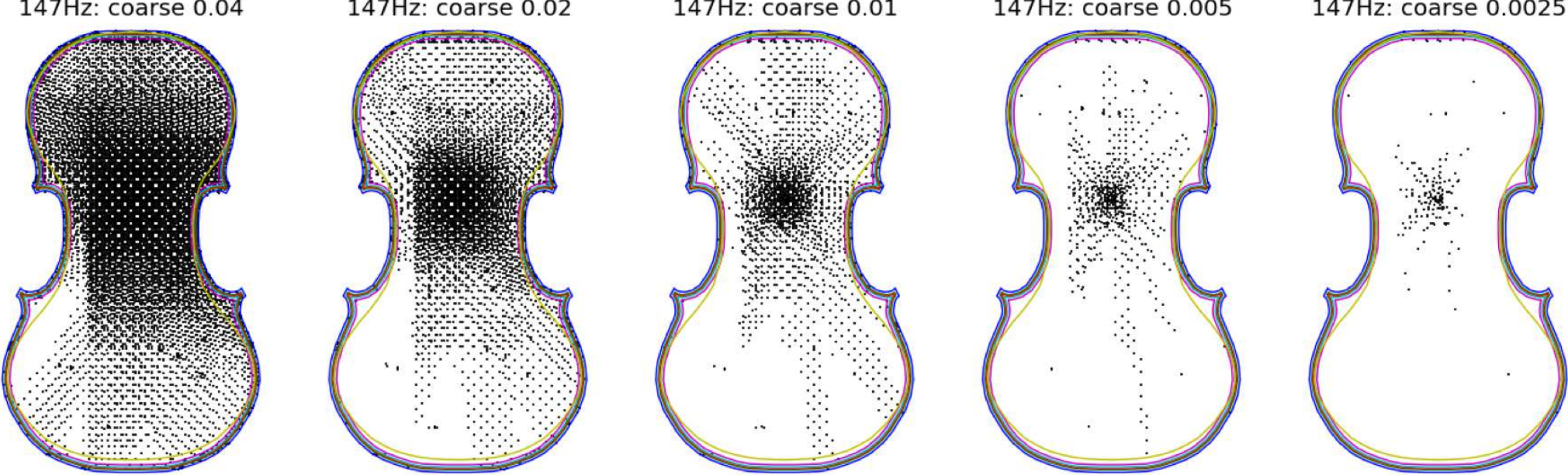}
\smallskip

\includegraphics[height=5.2in,angle=-90]{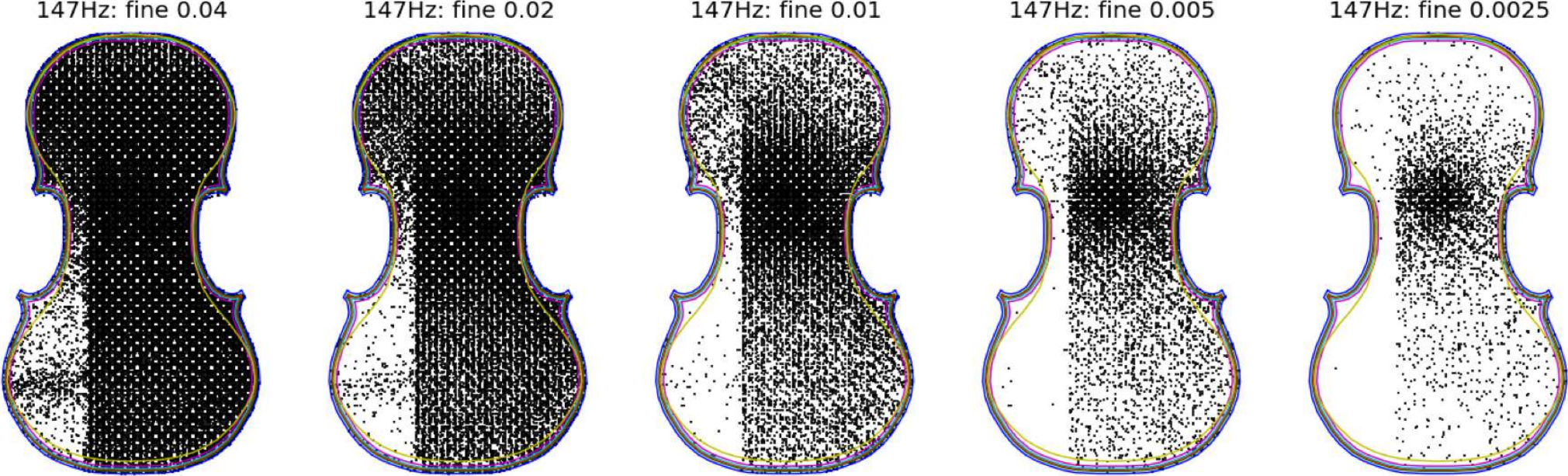}
\smallskip

\caption{Nodal points of the resonance wave at $f=147$Hz arising from
the six mode of vibrations of the plate of frequencies closest to $f$. From 
left to right, the notion of nodal point is 
defined by choosing $c_\Omega=0.04, 0.02, 0.01, 0.005$ and $0.0025$, 
respectively.}
\label{f7}
\end{figure}

\end{document}